\documentclass[10pt]{article}

\RequirePackage[amsthm,amsmath]{imsart}

\usepackage[mathscr]{eucal}
\usepackage{graphicx}
\usepackage{latexsym}
\usepackage{amssymb}
\usepackage{psfrag}
\usepackage{epsfig}
\usepackage{enumerate}
\DeclareMathAlphabet{\mathpzc}{OT1}{pzc}{m}{it}

\include{psbox}

\usepackage{amsfonts}
\usepackage{graphicx}
\usepackage{amsfonts}
\usepackage{color}
\usepackage[usenames,dvipsnames]{xcolor}

\usepackage[numbers]{natbib}

\begin{document}

	\newtheorem{proposition}{Proposition}[section]
	\newtheorem{theorem}[proposition]{Theorem}
	\newtheorem{corollary}[proposition]{Corollary}
	\newtheorem{lemma}[proposition]{Lemma}
	\newtheorem{conjecture}[proposition]{Conjecture}
	\newtheorem{question}[proposition]{Question}
	\newtheorem{definition}[proposition]{Definition}
	\newtheorem{comment}[proposition]{Comment}
	\newtheorem{algorithm}[proposition]{Algorithm}
	\newtheorem{assumption}[proposition]{Assumption}
	\newtheorem{condition}[proposition]{Condition}
	\numberwithin{equation}{section}
	\numberwithin{proposition}{section}

\newcommand{\skp}{\vspace{\baselineskip}}
\newcommand{\noi}{\noindent}
\newcommand{\osc}{\mbox{osc}}
\newcommand{\lfl}{\lfloor}
\newcommand{\rfl}{\rfloor}

\theoremstyle{remark}
\newtheorem{example}{\bf Example}[section]
\newtheorem{remark}{\bf Remark}[section]

\newcommand{\img}{\imath}
\newcommand{\iy}{\infty}
\newcommand{\eps}{\varepsilon}
\newcommand{\del}{\delta}
\newcommand{\Rk}{\mathbb{R}^k}
\newcommand{\RR}{\mathbb{R}}
\newcommand{\spa}{\vspace{.2in}}
\newcommand{\V}{\mathcal{V}}
\newcommand{\E}{\mathbb{E}}
\newcommand{\I}{\mathbb{I}}
\newcommand{\PP}{\mathbb{P}}
\newcommand{\sgn}{\mbox{sgn}}
\newcommand{\ti}{\tilde}

\newcommand{\QQ}{\mathbb{Q}}

\newcommand{\XX}{\mathbb{X}}
\newcommand{\XXz}{\mathbb{X}^0}

\newcommand{\lan}{\langle}
\newcommand{\ran}{\rangle}
\newcommand{\lf}{\lfloor}
\newcommand{\rf}{\rfloor}
\def\wh{\widehat}
\newcommand{\defn}{\stackrel{def}{=}}
\newcommand{\txb}{\tau^{\epsilon,x}_{B^c}}
\newcommand{\tyb}{\tau^{\epsilon,y}_{B^c}}
\newcommand{\tilxb}{\tilde{\tau}^\eps_1}
\newcommand{\pxeps}{\mathbb{P}_x^{\eps}}
\newcommand{\non}{\nonumber}
\newcommand{\dist}{\mbox{dist}}

\newcommand{\Om}{\mathnormal{\Omega}}
\newcommand{\om}{\omega}
\newcommand{\vph}{\varphi}
\newcommand{\Del}{\mathnormal{\Delta}}
\newcommand{\Gam}{\mathnormal{\Gamma}}
\newcommand{\Sig}{\mathnormal{\Sigma}}

\newcommand{\tilyb}{\tilde{\tau}^\eps_2}
\newcommand{\beq}{\begin{eqnarray*}}
\newcommand{\eeq}{\end{eqnarray*}}
\newcommand{\beqn}{\begin{eqnarray}}
\newcommand{\eeqn}{\end{eqnarray}}
\newcommand{\ink}{\rule{.5\baselineskip}{.55\baselineskip}}

\newcommand{\bt}{\begin{theorem}}
\newcommand{\et}{\end{theorem}}
\newcommand{\deps}{\Del_{\eps}}
\newcommand{\dbl}{\mathbf{d}_{\tiny{\mbox{BL}}}}

\newcommand{\be}{\begin{equation}}
\newcommand{\ee}{\end{equation}}
\newcommand{\ac}{\mbox{AC}}
\newcommand{\BB}{\mathbb{B}}
\newcommand{\VV}{\mathbb{V}}
\newcommand{\DD}{\mathbb{D}}
\newcommand{\KK}{\mathbb{K}}
\newcommand{\HH}{\mathbb{H}}
\newcommand{\TT}{\mathbb{T}}
\newcommand{\CC}{\mathbb{C}}
\newcommand{\ZZ}{\mathbb{Z}}
\newcommand{\SSS}{\mathbb{S}}
\newcommand{\EE}{\mathbb{E}}
\newcommand{\NN}{\mathbb{N}}

\newcommand{\clg}{\mathcal{G}}
\newcommand{\clb}{\mathcal{B}}
\newcommand{\cls}{\mathcal{S}}
\newcommand{\clc}{\mathcal{C}}
\newcommand{\clj}{\mathcal{J}}
\newcommand{\clm}{\mathcal{M}}
\newcommand{\clx}{\mathcal{X}}
\newcommand{\cld}{\mathcal{D}}
\newcommand{\cle}{\mathcal{E}}
\newcommand{\clv}{\mathcal{V}}
\newcommand{\clu}{\mathcal{U}}
\newcommand{\clr}{\mathcal{R}}
\newcommand{\clt}{\mathcal{T}}
\newcommand{\cll}{\mathcal{L}}
\newcommand{\clz}{\mathcal{Z}}

\newcommand{\cli}{\mathcal{I}}
\newcommand{\clp}{\mathcal{P}}
\newcommand{\cla}{\mathcal{A}}
\newcommand{\clf}{\mathcal{F}}
\newcommand{\clh}{\mathcal{H}}
\newcommand{\N}{\mathbb{N}}
\newcommand{\Q}{\mathbb{Q}}
\newcommand{\bfx}{{\boldsymbol{x}}}
\newcommand{\bfa}{{\boldsymbol{a}}}
\newcommand{\bfh}{{\boldsymbol{h}}}
\newcommand{\bfs}{{\boldsymbol{s}}}
\newcommand{\bfm}{{\boldsymbol{m}}}
\newcommand{\bff}{{\boldsymbol{f}}}
\newcommand{\bfb}{{\boldsymbol{b}}}
\newcommand{\bfw}{{\boldsymbol{w}}}
\newcommand{\bfz}{{\boldsymbol{z}}}
\newcommand{\bfu}{{\boldsymbol{u}}}
\newcommand{\bfell}{{\boldsymbol{\ell}}}
\newcommand{\bfn}{{\boldsymbol{n}}}
\newcommand{\bfd}{{\boldsymbol{d}}}
\newcommand{\bfbeta}{{\boldsymbol{\beta}}}
\newcommand{\bfzeta}{{\boldsymbol{\zeta}}}
\newcommand{\bfnu}{{\boldsymbol{\nu}}}
\newcommand{\bfvarphi}{{\boldsymbol{\varphi}}}

\newcommand{\curvz}{{\bf \mathpzc{z}}}
\newcommand{\curvx}{{\bf \mathpzc{x}}}
\newcommand{\curvi}{{\bf \mathpzc{i}}}
\newcommand{\curvs}{{\bf \mathpzc{s}}}
\newcommand{\blip}{\mathbb{B}_1}

\newcommand{\BM}{\mbox{BM}}

\newcommand{\tac}{\mbox{\scriptsize{AC}}}



\begin{frontmatter}
\title{Central Limit Results for  Jump-Diffusions with Mean Field Interaction and a Common Factor.}

 \runtitle{CLT for Weakly Interacting Particles}

\begin{aug}
\author{ Amarjit Budhiraja\thanks{Research  supported in part by the National Science Foundation (DMS-1004418, DMS-1016441, DMS-1305120) and the Army Research
 Office (W911NF-10-1-0158, W911NF- 14-1-0331)}, Elisabeti Kira\thanks{Research support in part by  CAPES and Fulbright Commission (Proc. 15856/12-7)} and Subhamay Saha\thanks{Research supported in part by the Indo-US VI-MSS postdoctoral fellowship}\\ \ \\
}
\end{aug}

\today

\skp

\begin{abstract}
A system of $N$ weakly interacting particles whose dynamics is given in terms of jump-diffusions with a common factor is considered.  The common factor is described through 
another jump-diffusion and the coefficients of the evolution equation for each particle depend, in addition to its own state value, on the empirical
measure of the states of the $N$ particles and the common factor.  A Central Limit Theorem, as $N \to \infty$, is established.  The limit law is described in terms of a certain Gaussian mixture.
An application to models in Mathematical Finance of self-excited correlated defaults  is described.

\noi {\bf AMS 2000 subject classifications:} 60F05; 60K35; 60H30; 60J70.

\noi {\bf Keywords:} Mean field interaction, common factor, weakly interacting jump-diffusions, propagation of chaos, central limit theorems, fluctuation limits, symmetric statistics, multiple Wiener integrals,
self-excited correlated defaults.
\end{abstract}

\end{frontmatter}

\section{Introduction}\label{introsec}

For $N \ge 1$, let $Z^{N,1}, \cdots Z^{N,N}$ be $\RR^d$ valued stochastic processes, representing trajectories of $N$ particles, which are described through stochastic differential equations (SDE) driven by mutually independent Brownian motions(BM) and Poisson random measures(PRM) such that the statistical distribution of $(Z^{N,1}, \ldots, Z^{N,N})$ is exchangeable.  The dependence
between the $N$ stochastic processes enters through the coefficients of the SDE which, for the $i$-th process, depend in addition
to the $i$-th state process, on a common stochastic process (common factor) and the empirical measure
$\mu^N_t = \frac{1}{N} \sum_{i=1}^N \delta_{Z^{N,i}_t}$.  The common factor is a $m$-dimensional stochastic process described
once more through a SDE driven by a BM and a PRM which are independent of the other noise processes.  Such stochastic systems are commonly referred to as {\em weakly interacting Markov processes} and have a long history.
Some of the classical works include McKean\cite{McK1, McK2},  Braun and Hepp \cite{BH}, Dawson \cite{Daw}, Tanaka \cite{Tan}, Oelschal\"{a}ger \cite{Oel}, Sznitman \cite{Sznit-a, Sznit}, Graham and M\'{e}l\'{e}ard \cite{GraMel}, Shiga and Tanaka \cite{ShiTan}, M\'{e}l\'{e}ard \cite{Mel}.  All of these papers treat the setting where the `common factor' is absent.  Most of this research activity
is centered around proving Law of Large Number results and Central Limit Theorems(CLT).  For example one can show (cf. \cite{Sznit-a, Oel}) that under suitable conditions, if the joint initial distributions of every set of  $k$-particles, for every $k$, converge to  product measures as $N \to \infty$ then the same is true for the joint distribution of the 
 stochastic processes(considered as path space valued random variables) as well.  Such a result, referred to as the {\em propagation of chaos} is one of the key first steps in the study of the fluctuation theory for such a system of interacting particles.

Systems with a common factor arise in many different areas.  In Mathematical Finance, they have been used to model 
correlations between default probabilities of multiple firms\cite{CMZ}.  In neuroscience modeling these arise as systematic noise 
in the external current input to a neuronal ensemble\cite{FTC}.  For particle approximation schemes for stochastic
partial differential equations (SPDE), the common factor corresponds to the underlying driving noise in the 
SPDE\cite{KuXi1, KuXi2}.
The goal of this work is to study a general family of  weakly interacting jump-diffusions with a common factor.
Our main objective is to establish a suitable Central Limit Theorem.  A key point here is that due to the presence of the common factor, the limit of $\frac{1}{N} \sum_{i=1}^N \delta_{Z^{N,i}}$ will in general be a random measure.  This in particular
means that the centering in the fluctuation theorem  will typically be random as well and one expects the limit law for such
fluctuations to be not Gaussian but rather a `Gaussian mixture'.  Our main result is Theorem \ref{thm:clt}
which provides a CLT under Conditions \ref{cond:main}, \ref{cond:clt1}, \ref{cond:clt2}, \ref{cond:clt3}, and \ref{cond:clt4}.  The summands in this CLT can be quite general
functionals of the trajectories of the particles with suitable integrability properties. The key idea is to first consider
a closely related collection of $N$ stochastic processes that, conditionally on a common factor, are independent and identically distributed.  By introducing a suitable Radon-Nikodym derivative one can evaluate the expectations associated with 
a perturbed form of the original
scaled and centered sum in terms of the conditionally i.i.d. collection.  The asymptotics of the latter quantity are easier
to analyze using, in particular, the classical limit theorems for symmetric statistics\cite{DynMan}.  
The perturbation arises due to the fact that in the original system the evolution of the common factor jump-diffusion depends on the empirical measure
of the states of the $N$-particles whereas in the conditionally i.i.d. construction the common factor evolution is determined by the large particle limit of the empirical measures.  Estimating the error introduced by this perturbation is one of the
key technical challenges in the proof.

In a setting where there is no common factor such central limit results have been obtained
in the classical works of Sznitman\cite{Sznit} and Shiga and Tanaka\cite{ShiTan}.  In this case the limit law is Gaussian and the
probability law of the actual $N$-particle system can be realized exactly through a simple absolutely continuous change of measure from the probability law of an i.i.d. system.   Another aspect that makes the analysis in the current work significantly more challenging is that unlike \cite{Sznit, ShiTan} 
 the dependence of the coefficients of the model  on the empirical measure in nonlinear.  

Central limit theorems for  systems of weakly interacting particles with a common factor have previously been
studied in \cite{KuXi2}.  This work is motivated by applications to particle system approximations to solutions of SPDE.
In addition to the fact that the form of the common factor in \cite{KuXi2} is quite different from that in our work, there are several differences between these two works.   The model considered in the current work allows for jumps in both
particle dynamics and the common factor dynamics neither of which are present in \cite{KuXi2}.  Also, in \cite{KuXi2} the fluctuation limit theorem is established for centered and scaled empirical measures considered as stochastic processes in the space of (modified) Schwartz distributions which in practice yields a functional central limit theorem for smooth functionals that depend on just the current state
of the particles.  In contrast, the current work allows for very general square integrable functionals that could possibly depend
on the whole trajectory of the particles. Thus, in particular, unlike \cite{KuXi2}, one can obtain from our work  limit theorems for statistics that depend on the particle states at multiple time instants.
In Section \ref{sec:pathspace} we sketch an argument that shows how one can recover convergence of modified Schwartz distribution valued stochastic processes from our main convergence result (Theorem \ref{thm:clt}).
A key difference in the argument here (from \cite{KuXi2}) is that we do not require unique solvability results for SPDE in order to characterize the limit.
More precisely, in \cite{KuXi2} the limit law is characterized through the solution of a certain SPDE and one of the key technical challenges is proving the wellposedness of the equation, whereas in the current work the description of the limit law  is given in terms of a certain mixture of Gaussian distributions (see \eqref{eq:limlaw}).  We note that in some respects the results in \cite{KuXi2} are more general in that
they allow for    infinite dimensional common factors 
and weighted empirical measures.  Our proofs rely on a Girsanov change of measure which requires the diffusion coefficients to satisfy a suitable non-degeneracy condition. Although the proofs in \cite{KuXi2}
are quite different and the form of state dependence allowed there is somewhat more general, it is interesting to note that the approach taken in \cite{KuXi2} also requires a non-degeneracy condition on the diffusion coefficient (see Condition (S4) in Section 4 of \cite{KuXi2}).

One of our motivations for the current study is to establish central limit results for models in Mathematical Finance of
self-exciting correlated defaults\cite{CMZ}.  In Section \ref{sec:finapp} we describe how Theorem \ref{thm:clt} 
yields such results.

The paper is organized as follows.  In Section \ref{prelim} we begin by introducing our model
of weakly interacting jump-diffusions with a common factor.  Next in Section \ref{sec:wellpos} we present a basic condition(Condition \ref{cond:main})
that will ensure pathwise existence and uniqueness of solutions to the SDE for the $N$-particle system and also for a related family of SDE describing a nonlinear Markov process.  These wellposedness results are given in Theorem \ref{thm:wellpos}
the proof of which is given in Section \ref{sec:infsyst}. The proofs are based on ideas and results from
\cite{KurPro,KotKur,KuXi2}.  
In Section \ref{sec:cltassu} we present the main Central Limit Theorem of this work.  Introducing conditions for the theorem requires some notation and thus we postpone some of them to later sections
(specifically Conditions \ref{cond:clt2}, \ref{cond:clt3} and \ref{cond:clt4} are introduced in Sections \ref{sec:estyny}, \ref{sec:7.1} and \ref{sec:secjn2} respectively).
  Section \ref{sec:infsyst} is devoted to the proof of Theorem \ref{thm:wellpos}.  In Section \ref{sec:symmstat} we recall the classical result of Dynkin and Mandelbaum\cite{DynMan} on limit laws of degenerate symmetric
statistics described in terms of multiple Wiener integrals.  Section \ref{sec:gir} introduces the Girsanov change of measure that is the key ingredient in our proofs. Section \ref{sec:estyny} enables the 
estimation of the error due to the perturbation described earlier in the Introduction and 
Section \ref{sec:proofmainth} contains the proof of Theorem \ref{thm:clt}. In Section \ref{sec:pathspace}, using Theorem \ref{thm:clt}, we sketch an argument for proving weak convergence of scaled and centered empirical measures as stochastic
processes with values in the dual of a suitable Nuclear space.
 Finally Section \ref{sec:finapp} discusses an application of Theorem \ref{thm:clt} to certain models in Mathematical Finance.

The following notations will be used.  
Fix $T < \infty$. All stochastic processes will be considered over the time horizon $[0,T]$.  We will use the notations $\{X_t\}$ and $\{X(t)\}$ interchangeably for stochastic processes.
Space of probability measures on a Polish space $\SSS$, equipped with the topology of weak convergence, will be denoted by $\clp(\SSS)$.  A convenient metric for this topology 
is the bounded-Lipschitz metric $\dbl$ defined as
$$
\dbl(\nu_1, \nu_2) = \sup_{f \in \blip} | \langle f , \nu_1 - \nu_2\rangle|, \nu_1, \nu_2 \in \clp(\SSS),$$
where  $\blip$ is the collection of  all Lipschitz functions $f$ that are bounded by $1$ and such that the corresponding Lipschitz constant is bounded by $1$ as well; and 
$\langle f, \mu\rangle = \int f d\mu$ for a signed measure $\mu$ on $\SSS$ and $\mu$-integrable $f: \SSS \to \RR$.
For a function $f: [0,T] \to \RR^k$, $\|f\|_{*,t} \doteq \sup_{0\le s \le t}\|f(s)\|$, $t \in [0,T]$.  Also, for $\mu_i:[0,T] \to \clp(\SSS)$, $i=1,2$,
$$\dbl(\mu_1, \mu_2)_{*,t} = \sup_{0\le s \le t}\dbl(\mu_1(s), \mu_2(s)).$$
Borel $\sigma$-field on a Polish space $\SSS$ will be denoted as $\clb(\SSS)$.  Space of functions that are right continuous with left limits (RCLL) from $[0, \infty)$ [resp. $[0,T]$] to $\SSS$
will be denoted as $\DD_{\SSS}[0,\infty)$ [resp. $\DD_{\SSS}[0,T]$] and equipped with the usual Skorohod topology.
Similarly $\CC_{\SSS}[0,\infty)$ [resp. $\CC_{\SSS}[0,T]$] will be the space of continuous functions
from $[0, \infty)$ [resp. $[0,T]$] to $\SSS$,  equipped with the local uniform [resp. uniform] topology.
For $x \in \DD_{\SSS}[0,T]$ and $t \in [0, T]$, $x_{[0,t]}$ will denote the element of $\DD_{\SSS}[0,t]$ defined
as $x_{[0,t]}(s) = x(s)$, $s \in [0,t]$.
Also given $x_{[0,t]} \in \DD_{\SSS}[0,t]$, $x_{[0,t]}(s)$ will be written as $x_s$.  Similar notation will be used for stochastic processes.

For a bounded function $f$ from $\SSS$ to $\RR$,  $\|f\|_{\infty} = \sup_{x \in \SSS}|f(x)|$.
Probability law of a $\SSS$ valued random variable $\eta$ will be denoted as $\cll(\eta)$ and its conditional distribution (a $\clp(\SSS)$ valued random variable) given 
a sub-$\sigma$ field $\clg$ will be denoted as $\cll(\eta \mid \clg)$. 
Convergence of a sequence $\{X_n\}$ of $\SSS$ valued random variables in distribution to $X$ will be written as $X_n \Rightarrow X$.
For a $\sigma$-finite measure $\nu$ on a Polish space $\SSS$, $L^2_{\RR^k}(\SSS, \nu)$ will denote the Hilbert space of $\nu$-square integrable functions from $\SSS$ to $\RR^k$.  When $k=1$, we will merely write $L^2(\SSS, \nu)$.
The norm in this Hilbert space will be denoted as $\|\cdot\|_{L^2(\SSS, \nu)}$.
We will usually denote by $\kappa, \kappa_1, \kappa_2, \cdots$, the constants that appear in various estimates within a proof.  The values of these constants may change from one proof to another.

\section{Main results}\label{prelim}
Let
$$\XX_t = [0,t] \times \RR^d \times \RR_+,\;\; \XXz_t = [0,t] \times \RR^m \times \RR_+, \; t \in [0,T].$$
For $k \in \N$, let $\clc_k$ and  $\cld_k$ denote $\CC_{\RR^k}[0,T]$ and $\DD_{\RR^k}[0,T]$ respectively.
Let $\clm_d$ [resp. $\clm_m$]  be the space of $\sigma$-finite measures on $\XX_T$ [resp. $\XX_T^0$] with the topology of vague convergence.

For fixed $N \ge 1$, consider the system of equations for the $\RR^d$ valued RCLL stochastic processes $Z^{N,i}$, $i = 1, \ldots N$ and the $\RR^m$ valued 
RCLL process $U^N$ given on a filtered probability space $(\Om, \clf, \PP, \{\clf_t\})$:
\begin{align}
	Z^{N,i}_t &= Z^{N,i}_0 + \int_0^t b(Z^{N,i}_s, U^N_s, \mu^N_s) ds + \int_0^t \sigma(Z^{N,i}_s, U^N_s, \mu^N_s) dB^i_s\nonumber\\
&\quad + \int_{\XX_t}\psi_d(Z^{N,i}_{s-}, U^N_{s-}, \mu^N_{s-},u,h) d\bfn^i\label{eq:nparteq}\\
U^N_t &= U_0 + \int_0^t b_0(U^N_s, \mu^N_s) ds + \int_0^t \sigma_0(U^N_s, \mu^N_s) dB^0_s\nonumber\\
&\quad + \int_{\XX^0_t}\psi_{d_0}(U^N_{s-}, \mu^N_{s-},u,h) d\bfn^0,\label{eq:comminp}
\end{align}
Here $B^i$, $i\in \NN$ are $r$-dimensional Brownian motions(BM); $B^0$ is a $m$ dimensional BM; $\bfn^i$, $i\in \NN$  are Poisson random measures (PRM)
with intensity measure $\bfnu = \lambda_T \otimes \gamma \otimes \lambda_{\infty}$ on $\XX_T$, where $\lambda_T$ [resp. $\lambda_{\infty}$] is the Lebesgue measure on $[0,T]$ [resp.
$[0,\infty)$] and $\gamma$ is a finite measure on $\RR^d$; $\bfn^0$ is a PRM with intensity measure $\bfnu^0 = \lambda_T \otimes \gamma^0 \otimes \lambda_{\infty}$ on $\XXz_T$, where  $\gamma^0$ is a finite measure on $\RR^m$.  All these processes are mutually independent and they have independent increments with respect to the filtration $\{\clf_t\}$.
Also $\mu^N_s = \frac{1}{N}\sum_{i=1}^N \delta_{Z^{N,i}_s}$ and $\psi_d, \psi_{d_0}$ are maps defined as follows: For 
$(x,y,\theta,u,h,k) \in \RR^d\times \RR^m\times \clp(\RR^d)\times
\RR_+ \times \RR^d \times \RR^m$
$$
\psi_d(x,y,\theta,u,h) = h1_{[0, d(x,y,\theta,h)]}(u), \; \psi_{d_0}(y,\theta,u,k) = k1_{[0, d_0(y,\theta,k)]}(u),$$
where $d$ and $d_0$ are nonnegative maps on $\RR^{d+m}\times \clp(\RR^d) \times \RR^d$ and $\RR^{m}\times \clp(\RR^d) \times \RR^m$ respectively.
Roughly speaking, given $(Z^{N,i}_{t-}, U^N_{t-}, \mu^N_{t-}) = (x,y,\theta)$, the jump for $Z^{N,i}$ at instant $t$
occurs at rate $\int_{\RR^d} d(x,y,\theta,h) \gamma(dh)$ and the jump distribution is given as
$c\cdot d(x,y,\theta,h) \gamma(dh)$ where $c$ is the normalization constant. Jumps of $U^N$ are described in an analogous manner.

We assume that $\{Z_0^{N,i}\}_{i=1}^N$ are i.i.d. with common distribution $\mu_0$ and $U_0$ is independent of $\{Z_0^{N,i}\}_{i=1}^N$ and has probability distribution $\rho_0$.  Also, $\{Z_0^{N,i}\}_{i=1}^N$ and $U_0$ are 
$\clf_0$ measurable.

Conditions on the various coefficients will be introduced shortly.  
Along with the $N$-particle equations \eqref{eq:nparteq}-\eqref{eq:comminp} we will also consider  a related infinite system of equations for 
$\RR^d \times \RR^m$  valued RCLL stochastic processes $(X^i,Y)$, $i \in \NN$ given on   $(\Om, \clf, \PP, \{\clf_t\})$.
\begin{align}
	X^i_t &= X^i_0 + \int_0^t b(X^i_s, Y_s, \mu_s) ds + \int_0^t \sigma(X^i_s, Y_s, \mu_s) dB^i_s\nonumber\\
&\quad + \int_{\XX_t}\psi_d(X^i_{s-}, Y_{s-}, \mu_{s-},u,h) d\bfn^i\label{eq:mfield1}\\
Y_t &= Y_0 + \int_0^t b_0(Y_s, \mu_s) ds + \int_0^t \sigma_0(Y_s, \mu_s) dB^0_s\nonumber\\
&\quad + \int_{\XX^0_t}\psi_{d_0}(Y_{s-}, \mu_{s-},u,h) d\bfn^0,\label{eq:mfield2}
\end{align}
Here 
$\mu_t = \lim_{k\to \infty} \frac{1}{k} \sum_{i=1}^k \delta_{X^i_t}$, where the limit is a.s. in $\clp(\RR^d)$.
%
%
%
As for the $N$-particle system, we assume that $\{X_0^{i}\}_{i\in \mathbb{N}}$ are i.i.d. with common distribution $\mu_0$ and $Y_0$ is independent of $X_0 \equiv \{X_0^{i}\}_{i\in \mathbb{N}}$ and has probability distribution $\rho_0$.  Also, $\{X_0^{i}\}_{i\in \mathbb{N}}$ and $Y_0$ are $\clf_0$ measurable.

\subsection{Well-posedness.}\label{sec:wellpos}
We now give conditions on the coefficient functions under which the systems of equations \eqref{eq:nparteq}-\eqref{eq:comminp} and
\eqref{eq:mfield1}-\eqref{eq:mfield2} have unique pathwise solutions.  A pathwise solution of \eqref{eq:mfield1}-\eqref{eq:mfield2} is a collection of
RCLL processes $(X^i,Y)$, $i \ge 1$,  with values in $\RR^d\times \RR^m$  such that: (a) $Y$ is $\{\clg^0_t\}$ adapted, where
$\clg^0_s = \sigma \{ Y_0, B^0_r, \bfn^0([0,r]\times A), r \le s, A \in \clb(\RR^m \times \RR_+ )\}$;
(b) $X$ is $\{\clf_t\}$ adapted where $X = (X^i)_{i \in \NN}$;
(c) stochastic integrals on the right sides of \eqref{eq:mfield1}-\eqref{eq:mfield2} are well defined; (d) Equations \eqref{eq:mfield1}-\eqref{eq:mfield2} hold a.s.
Uniqueness of pathwise solutions says that if $(X,Y)$ and $(X', Y')$ are two such solutions with $(X_0,Y_0) = (X'_0,Y'_0)$ then they must be indistinguishable.
Existence and uniqueness of solutions to \eqref{eq:nparteq}-\eqref{eq:comminp} are defined in a similar manner.  In particular, in this case (a) and (b) are replaced by the requirement that $(Z^{N,i}, U^N)_{i=1}^N$ are $\{\clf_t\}$ adapted.

We now introduce conditions on the coefficients that will ensure existence and uniqueness of solutions.
\begin{condition}
	\label{cond:main}
	There exist   $\epsilon, K \in (0, \infty)$ such that
	\begin{description}
		\item{(a)} For all $z=(x,y) \in \RR^d \times \RR^m$, $\nu \in \clp(\RR^d)$, $(h,k)\in \RR^d \times \RR^m$,
		$$\epsilon \le d(z,\nu,h) \le K, \; 0\le d_0(y,\nu,k)< K, \; \int_{\RR^d} \|h\|^2 \gamma(dh) \le K^2, \; \int_{\RR^m} \|k\|^2 \gamma^0(dk) \le K^2,$$
		and
		$$ \max\left\{ \|\sigma(z,\nu)\|,\, \|\sigma_0(y,\nu)\|,\,  \|b(z,\nu)\|,\, \|b_0(y,\nu)\|\right\} \le K.$$
		\item{(b)} For all $z=(x,y), z'=(x',y') \in \RR^d \times \RR^m$, $\nu,\nu' \in \clp(\RR^d)$ the functions $\sigma, \sigma_0, b, b_0$ satisfy
		\begin{align*}
		\|\sigma(z,\nu) - \sigma(z',\nu')\| + \|\sigma_0(y,\nu) - \sigma_0(y',\nu')\|  &\le K (\|z-z'\| + \dbl(\nu,\nu'))\\
		\|b(z,\nu) - b(z',\nu')\| +  \|b_0(y,\nu) - b_0(y',\nu')\| &\le K (\|z-z'\| + \dbl(\nu,\nu'))
		\end{align*}
		and the functions $d, d_0$ satisfy
		\begin{align*}
		\int_{\RR^d}  \|h\|^2 \|d(z,\nu, h) - d(z',\nu', h)\| \gamma(dh) &\le K (\|z-z'\| + \dbl(\nu,\nu'))\\
		 \int_{\RR^m}  \|k\|^2 \|d_0(y,\nu, k) - d_0(y',\nu', k)\| \gamma^0(dk)  &\le K (\|y-y'\| + \dbl(\nu,\nu'))
	\end{align*}
	\end{description}
	
\end{condition}
Under the above condition we can establish the following wellposedness result.  
\begin{theorem}
\label{thm:wellpos}
 Suppose that 
\begin{equation}
	\label{eq:eq734}
	\int \|x\|^2 \mu_0(dx) + \int \|y\|^2 \rho_0(dy) < \infty
	\end{equation}
  and  Condition \ref{cond:main} holds. Then: \\

\noindent (a) the system of equations \eqref{eq:mfield1}-\eqref{eq:mfield2} has a unique pathwise solution. \\ 

\noindent (b)  the system of equations \eqref{eq:nparteq}-\eqref{eq:comminp} has a unique pathwise solution.  
\end{theorem}

Proof of the theorem  is given in Section \ref{sec:infsyst}.  
\begin{remark}
	\label{rem:rem415}
	(i) We note that the unique pathwise solvability in (a) implies that there is a measurable map $\clu: \RR^m \times \clc_m \times \clm_m \to \cld_m$ such that
	the solution $Y$ of \eqref{eq:mfield2} is given as $Y = \clu(Y_0, B^0, \bfn^0)$.\\
(ii) Recall that $\clg^0_s = \sigma \{Y_0, B^0_r, \bfn^0([0,r]\times A), r \le s, A \in \clb(\RR^m \times \RR_+ )\}$, $s \in [0, T]$. Let $\clg^0 = \clg^0_T$.
Then exactly along the lines of Theorem 2.3 of \cite{KuXi1} it follows that if $(\{X^i\}, Y)$ is a solution of \eqref{eq:mfield1}-\eqref{eq:mfield2} then
\begin{equation}
	\mu_t = \cll(X^i(t) \mid \clg^0) = \cll(X^i(t) \mid \clg^0_t), \; t \in [0,T], \; i \in \NN .
\end{equation}
  In particular, there is a measurable map $\Pi: \mathbb{R}^m\times \clc_m \times \clm_m \to \mathbb{D}_{\clp(\RR^d)}[0,T]$ such that
$\Pi(Y_0,B^0,\bfn^0) = \mu$ a.s. 
\end{remark}
\subsection{Central Limit Theorem.}
\label{sec:cltassu}
The main result of this work establishes a CLT for $\frac{1}{N}\sum_{i=1}^N \delta_{Z^{N,i}}$.  
For that, we will make  additional assumptions on the coefficients.
\begin{condition}
	\label{cond:clt1}
	For some $p \in \NN$, $r= d+p$ and for all $(x,y, \nu, k) \in \RR^d\times \RR^m \times \clp(\RR^d) \times \RR^m$, $\sigma(x,y,\nu) = [I_{d\times d}, \; \tilde \sigma(x)]$,  $d_0(y,\nu,k) = d_0(k)$
	where $I_{d\times d}$ is the $d\times d$ identity matrix and $\tilde \sigma(x)$ is a $d\times p$ matrix.
\end{condition}
Note that $\sigma_0$ is allowed to depend on $(y,\nu)$. 
\begin{remark}
	\label{rem2.3mod}
	It is easily seen that if $b$ is of the form $b = (\tilde b, 0)'$ where for some $q < d$, $\tilde b$ is a $\RR^q$ valued function then one can relax the assumption on $\sigma$ by allowing it to be of the form
	$$\sigma(x, y, \nu) = \left( \begin{array}{cc}
	I_{q\times q} & \tilde \sigma_{12}(x)  \\
	\tilde \sigma_{21}(x) & \tilde \sigma_{22}(x)  \\
	 \end{array} \right).$$
\end{remark}

We will need additional smoothness assumptions on the coefficients $b, d, b_0$ and $\sigma_0$ (Conditions \ref{cond:clt2}, \ref{cond:clt3} and \ref{cond:clt4}) however stating them requires  some notation
which we prefer to introduce in later sections. As argued in Section \ref{sec:remoncond}, these conditions are satisfied quite generally.  Below is the main result of this work.
We begin by introducing the following canonical spaces and stochastic processes.
Let
$$\Omega_d = \clc_r \times \clm_d \times \cld_d, \;\; \Omega_m = \clc_m \times \clm_m \times \cld_m,$$
Recall from \eqref{eq:mfield1}-\eqref{eq:mfield2} the processes $(B^i, \bfn^i)_{i \in \NN_0}$ and the pathwise solution $(\{X^i\}_{i \in \NN},Y)$.
Define for $N \in \N$ the probability measure $P^N$ on $\bar{\Om}_{N} = \Om_{m} \times \Om_d^N$ as 
$$P^N = \cll\left (  (B^0, \bfn^0, Y), (B^1, \bfn^1, X^1), \ldots ,  (B^N, \bfn^N, X^N)\right)$$
Note that $P^N$ can be disintegrated as
\begin{equation}\label{eq:transkern}
	P^N(d\om_0\, d\om_1 \cdots d\om_N) = \alpha(\om_0,\, d\om_1) \cdots \alpha(\om_0, d\om_N) P_0(d\om_0),
\end{equation}
where $P_0 = \mathcal{L}(B^0, \bfn^0, Y)$.
For  $\bar \om = (\om_0, \om_1, \ldots , \om_N) \in \bar \Om_N$,  $V^i(\bar \om) = \om_i$, $i =0,1, \ldots , N$ and abusing notation,
\begin{equation}\label{eq:canon1} V^i = (B^i, \bfn^i, X^i), \; i = 1, \ldots , N, \;\; V^0 = (B^0, \bfn^0, Y).\end{equation}
Also define the canonical process $V_* = (B_*, \bfn_*, X_*)$  on $\Om_d$  as
\begin{align}
\label{eq:canon2}V_*(w) &= (B_*(w), \bfn_*(w), X_*(w)) = (w_{1},w_{2},w_{3});\; w = (w_{1},w_{2},w_{3}) \in \Om_d.
\end{align}
We denote by $\cla$ the collection of all measurable maps  $\varphi: \cld_d \to \mathbb{R}$ such that
$\varphi(X_*) \in L^2(\Om_d, \alpha(\om_0, \cdot))$ for $P_0$ a.e. $\om_0 \in \Om_m$.
For $\varphi \in \cla$ and $\om_0 \in \Om_m$, let 
\begin{equation}\label{eq:7.a}
m_{\varphi}(\om_0) = \int_{\Om_d} \varphi(X_*(\om_1)) \alpha(\om_0, d\om_1), \; \Phi_{\om_0} = \varphi(X_*) - m_{\varphi}(\om_0).
\end{equation}
Let for $\om_0 \in \Om_m$ and $\varphi \in \cla$, $\sigma^{\varphi}_{\om_0} \in \RR_+$ be defined through \eqref{eq:eq1209}. Denote by $\pi_{\om_0}^{\varphi}$ the normal distribution with mean $0$ and standard deviation $\sigma_{\om_0}^{\varphi}$.  Let
$\pi^{\varphi} \in \clp(\RR)$ be defined as 
\begin{equation}
	\label{eq:limlaw}\pi^{\varphi} = \int_{\Om_m} \pi^{\varphi}_{\om_0} P_0(d\om_0).
	\end{equation}  
	Finally with $\{Z^{N,i}\}_{i=1}^N$ as defined in \eqref{eq:nparteq}
and $\varphi \in \cla$, let 
$$\bar\clv_N^{\varphi} = \sqrt{N} \left (\frac{1}{N} \sum_{j=1}^N \varphi(Z^{N,j}) - m_{\varphi}(\bar V^0)\right),$$
where $\bar V^0 = (B^0, \bfn^0, \clu(U_0, B^0, \bfn^0))$ and $\clu$ is as introduced below Theorem \ref{thm:wellpos}.
Denote by $\pi_N^{\varphi} \in \clp(\RR)$ the probability distribution of $\bar\clv_N^{\varphi}$.
The following is the main result of this work.
\begin{theorem}
	\label{thm:clt}
	Suppose that Conditions \ref{cond:main}, \ref{cond:clt1}, \ref{cond:clt2}, \ref{cond:clt3} and \ref{cond:clt4} hold. Then, for all $\varphi \in \cla$, $\pi_N^{\varphi}$ converges weakly to $\pi^{\varphi}$ as $N \to \infty$.
\end{theorem}
Rest of the paper is organized as follows. In Section \ref{sec:infsyst} we present the proof of the wellposedness result (Theorem \ref{thm:wellpos}).
Section \ref{sec:symmstat} recalls some classical results of Dynkin and Mandelbaum\cite{DynMan} on limits of symmetric statistics.  In Section \ref{sec:gir} we introduce the Girsanov change of measure that plays a key role
in proofs and Section \ref{sec:estyny} gives some moment bounds that will be frequently appealed to in our proofs. Section 
\ref{sec:proofmainth} contains the proof of our main result (Theorem \ref{thm:clt}). In Section \ref{sec:pathspace} we discuss how  Theorem \ref{thm:clt}  can be used to prove central limit theorems for centered and scaled empirical measures. Finally Section \ref{sec:finapp} considers an application of our results to certain models in mathematical finance.

\section{Proof of Theorem \ref{thm:wellpos}.}
\label{sec:infsyst}
Proof of the theorem follows along the lines of \cite{KuXi1}, we sketch the argument for the
first statement in Theorem \ref{thm:wellpos} and omit the proof of the second statement.
Namely, we show now that if $\{X^i_0, i \in \NN\}$ and $Y$ are as defined below \eqref{eq:mfield2}; \eqref{eq:eq734} holds; and Condition \ref{cond:main} is satisfied, then the systems of equations \eqref{eq:mfield1}-\eqref{eq:mfield2} has a unique pathwise solution.
We first argue pathwise uniqueness.  Suppose that
$R = \{R^i = (X^i,Y), i \in \NN\}$ and $\tilde R = \{\tilde R^i = (\tilde X^i,\tilde Y), i \in \NN\}$ are two solutions of \eqref{eq:mfield1}-\eqref{eq:mfield2}
with $R_0 = \tilde R_0$.  Then using Condition \ref{cond:main} and standard maximal inequalities, for $t \in [0,T]$,
\begin{align*}
\EE\left \|\int_0^\cdot [\sigma(R^i_s,  \mu_s) - \sigma(\tilde R^i_s,  \tilde \mu_s)]  dB^i_s\right\|_{*,t}
&\le \kappa_1 K \EE \left [ \int_0^t (\|R^i -\tilde R^i\|_{*,s} + \dbl(\mu_{\cdot}, \tilde \mu_{\cdot})_{*,s})^2 ds\right]^{1/2}\\
& \le \kappa_1 K \sqrt{t} \EE(\|R^i -\tilde R^i\|_{*,t} + \dbl(\mu_{\cdot}, \tilde \mu_{\cdot})_{*,t}).
\end{align*}
Here, $\tilde \mu_t = \lim_{k\to \infty} \frac{1}{k} \sum_{i=1}^k \delta_{\tilde X^i_t}$ and $\kappa_1$ is a global constant.
Similarly, 
\begin{align*}
	&\EE\left\|\int_{[0,\cdot] \times \RR^d} [\psi_d(R^i_{s-}, \mu_{s-},u,h) - \psi_d(\tilde R^i_{s-},  \tilde \mu_{s-},u,h)] d\bfn^i\right\|_{*,t}\\
	&\quad\le \EE \int_{\XX_{t}} \|\psi_d(R^i_{s-}, \mu_{s-},u,h) - \psi_d(\tilde R^i_{s-},  \tilde \mu_{s-},u,h)\| d\bfn^i\\
	&\quad\le \EE \int_{[0,t]\times \RR^d} \|h\| |d(R^i_{s},\mu_{s},h) - d(\tilde R^i_{s},  \tilde \mu_{s},h)| \gamma(dh) ds\\
	&\quad\le \kappa_2 \int_{[0,t]} \EE(\|R^i -\tilde R^i\|_{*,s} + \dbl(\mu_{\cdot}, \tilde \mu_{\cdot})_{*,s}) ds,
\end{align*}
where the last inequality uses Condition \ref{cond:main}(b).
One has analogous estimates for terms involving $\sigma_0$, $d_0$, $b$ and $b_0$.
Also by Fatou's lemma,
\begin{align*}
	\EE \dbl(\mu_{\cdot}, \tilde \mu_{\cdot})_{*,s} =
	\EE \sup_{0 \le u \le s} \sup_{f \in \blip} | \langle f, \mu_u - \tilde \mu_u \rangle|
	 \le \liminf_{k\to \infty} \frac{1}{k} \sum_{i=1}^k \EE\|X^i - \tilde X^i\|_{*,s}
	\le \sup_i \EE\|X^i - \tilde X^i\|_{*,s}.
\end{align*}
Letting 
$$a_t = \sup_i \EE\|X^i - \tilde X^i\|_{*,t} + \EE\|Y - \tilde Y\|_{*,t}, \; t \in [0,T]$$
we then have from the above estimates that for some $\kappa_3 \in (0, \infty)$
$$a_t \le \kappa_3( \int_0^t a_s ds + \sqrt{t} a_t), \; t \in [0,T]$$
Taking $t$ sufficiently small we see now that $a_s =0$ for all $s \in [0,t]$.  A recursive argument then shows that $a_s = 0$ for all $s \in [0,T]$.  This completes the proof of uniqueness.

Next we prove existence of solutions.  We will use  ideas and results from \cite{KuXi1} (specifically 
Lemma 2.1 and  Theorem 2.2 therein).  
Define for
$t \in [0,T]$ and $n \ge 1$, $D^n(t) = \frac{\lfl nt \rfl}{n}$, 
$B^{n,i}_t = B^i_{\frac{\lfl nt \rfl}{n}}$, 
 $B^{n,0}_t = B^0_{\frac{\lfl nt \rfl}{n}}$,
 $\bfn^{n,i}(A_1 \times [0,t]) = \bfn^i(A_1 \times [0, \frac{\lfl nt \rfl}{n}])$, $\bfn^{n,0}(A_2 \times [0,t]) = \bfn^0(A_2 \times [0, \frac{\lfl nt \rfl}{n}])$,
 $A_1 \in \clb(\RR^d\times [0,T])$, $A_2 \in \clb(\RR^m\times [0,T])$.
Let $R^n \doteq (R^{n,i} = (X^{n,i}, Y^n), i \in \NN)$ be the solution of \eqref{eq:mfield1}-\eqref{eq:mfield2} with $dt$, $(B^i, B^0, \bfn^i, \bfn^0)$ and $\mu_t$ replaced by
$dD^n(t)$,  $(B^{n,i}, B^{n,0}, \bfn^{n,i}, \bfn^{n,0})$ and
$\mu^n_t = \lim_{k\to \infty} \frac{1}{k} \sum_{i=1}^k \delta_{X^{n,i}_t}$, respectively.  Note that the solution is determined recursively over intervals of length $1/n$ and $\mu^n_t$ is well defined for every $t \in [0,T]$ since  
$\lim_{k\to \infty} \frac{1}{k} \sum_{i=1}^k \delta_{X^{n,i}_t}$ exists a.s. from the exchangeability of $\{X^{n,i}_t, i \in \NN\}$ which in turn is  a consequence of
the exchangeability of $\{X^{i}_0, i \in \NN\}$.
Using the boundedness of the coefficients it is straightforward to check that
$$
\E \left ( \|R^{n,i}_{t + r} - R^{n,i}_t\| \mid \clf_t\right) \le \kappa_4 \frac{\lfl n(t+r)\rfl - \lfl nt \rfl}{n}, \; t \in [0,T-r], r \ge 0, i \in \NN,$$
where $\kappa_4$ is a constant independent of $n, i, t, r$.
It then follows that for each $i \in \NN$, $\{R^{n,i}\}_{n \in \NN}$ is tight in $\DD_{\RR^d \times \RR^m}[0, T]$.
This proves tightness of the sequence $\{R^{n}\}_{n \in \NN}$ in $\left(\DD_{\RR^d \times \RR^m}[0, T]\right)^{\otimes \infty}$.
A similar estimate as in the above display shows that for $i, j \in \NN$, $\{R^{n,i}+ R^{n,j}\}_{n \in \NN}$ is tight in $\DD_{\RR^d \times \RR^m}[0, T]$.
Thus we have that $\{R^{n}\}_{n \in \NN}$ is tight in $\DD_{(\RR^d \times \RR^m)^{\otimes \infty}}[0, T]$ (see for example \cite{EthKur}, Problems 3.11.22 and 3.11.23).
Let $\bar R \doteq \{\bar R^i = (\bar X^i, \bar Y)\}_{i \in \NN}$ denote a sub-sequential weak limit point.  Then $\{\bar X^i\}$ is exchangeable as well and so 
$\bar \mu_t = \lim_{k \to \infty} \frac{1}{k} \sum_{i=1}^k \delta_{\bar X^i_t}$ is well defined where the limit exists a.s.
From  Lemma 2.1 in \cite{KuXi1} (see also \cite{KotKur}) it now follows that (along the chosen subsequence)
$(R^n, \mu^n)$ converges in distribution to $(\bar R, \bar \mu)$, in $\DD_{(\RR^d \times \RR^m)^{\otimes \infty}\times \clp(\RR^d)}[0, T]$.

We note that $\psi_d$ regarded as a map from $\RR^d \times \RR^m \times \clp(\RR^d)$ to $L^2_{\RR^d}(\RR_+ \times \RR^d, \lambda_{\infty}\otimes \gamma)$
is a continuous map.  Indeed for $z = (x,y)$, $z'=(x',y') \in \RR^d \times \RR^m$ and $\nu, \nu' \in \clp(\RR^d)$
\begin{align*}
	\int_{\RR_+ \times \RR^d} \|\psi_d(z,\nu, u,h) - \psi_d(z',\nu', u,h)\|^2 du \gamma(dh) &=
	\int_{\RR^d} \|h\|^2 |d(z,\nu, h) - d(z',\nu', h)| \gamma(dh)\\
	&\le K (\|z-z'\| + \dbl(\nu, \nu'))
\end{align*}
where the last inequality is from Condition \ref{cond:main}.  Similarly $\psi_{d_0}$ is a continuous map from
$\RR^m  \times \clp(\RR^d)$ to $L^2_{\RR^m}(\RR_+ \times \RR^m, \lambda_{\infty}\otimes \gamma^0)$.
Fix $p \in \NN$, $\varphi_1, \cdots \varphi_p \in L^2_{\RR^d}(\RR_+ \times \RR^d, \lambda_{\infty}\otimes \gamma)$
and $\tilde \varphi_1, \cdots \tilde \varphi_p \in L^2_{\RR^m}(\RR_+ \times \RR^m, \lambda_{\infty}\otimes \gamma^0)$.
Let $\cli^{n,i}_{\varphi_j}(t) = \int_{\XX_t} \varphi_j(u,h) d\bfn^{n,i}$, $\cli^{n,i}_{\tilde \varphi_j}(t) = \int_{\XX^0_t} \tilde \varphi_j(u,k) d\bfn^{n,0}$, $j = 1, \ldots p$,
$t \in [0,T]$.  Fix $\ell \in \NN$.  Consider the vector of processes consisting of
$\sigma(X^{n,i}_{\cdot}, Y^n_{\cdot}, \mu^n_{\cdot})$, $b(X^{n,i}_{\cdot}, Y^n_{\cdot}, \mu^n_{\cdot})$, $\sigma_0(Y^n_{\cdot}, \mu^n_{\cdot})$,
$b_0(Y^n_{\cdot}, \mu^n_{\cdot})$, $B^{n,i}_{\cdot}$, $B^{n,0}_{\cdot}$, $\cli^{n,i}_{\varphi_j}$, $\cli^{n,i}_{\tilde \varphi_j}$,
$\psi_d(X^{n,i}_{\cdot}, Y^n_{\cdot}, {\bf \cdot})$,
$\psi_{d_0}(Y^n_{\cdot}, {\bf \cdot})$, $i \le \ell$, $j \le p$.
Then by the continuity of $b, b_0, \sigma, \sigma_0$ and the continuity property of $\psi_d$, $\psi_{d_0}$ noted above
this vector of processes converges in distribution in $\DD_{E}[0, T]$ to the vector of processes obtained by replacing
$(X^{n,i}, Y^n, \mu^n, B^{n,i}, B^{n,0}, \bfn^{n,i}, \bfn^{n,0})$ with
$(\bar X^{i}, \bar Y, \bar \mu, \bar B^{i}, \bar B^{0}, \bar \bfn^{i}, \bar \bfn^{0})$.  Here
$E = \RR^k \times L^2_{\RR^d}(\RR_+ \times \RR^d, \lambda_{\infty}\times \gamma) \times L^2_{\RR^m}(\RR_+ \times \RR^m, \lambda_{\infty}\times \gamma^0)$ for a suitable value
of $k$.  From Theorem 4.2 of \cite{KurPro} it now follows that  $(\bar X^i, \bar Y)$ is a solution of \eqref{eq:mfield1}-\eqref{eq:mfield2}  with 
$( B^{i},  B^{0},  \bfn^{i},  \bfn^{0})$ replaced with
$(\bar B^{i}, \bar B^{0}, \bar \bfn^{i}, \bar \bfn^{0})$ proving the existence of a weak solution of \eqref{eq:mfield1}-\eqref{eq:mfield2}.  From pathwise uniqueness
established earlier it now follows that there exists a strong solution of \eqref{eq:mfield1}-\eqref{eq:mfield2}.  
Exactly along the lines of the proof of Theorem 2.3 of \cite{KuXi1} it follows that $\{\mu_t\}$ is
$\{\clg^0_t\}$ adapted.  Also, using Condition \ref{cond:main}, if $(Y, \mu)$ and $(\tilde Y, \mu)$ solve \eqref{eq:mfield2}
then $Y$ and $\tilde Y$ are indistinguishable.  From this and the classical Yamada-Watanabe argument (cf. \cite{IkWa}, Theorem IV.1.1) it follows that $\{Y_t\}$ is $\{\clg^0_t\}$ adapted as well.
This completes the proof of pathwise existence and uniqueness of solutions. \qed


\section{Asymptotics of Symmetric Statistics.}
\label{sec:symmstat}
The proof of the central limit theorem crucially relies on certain classical results from \cite{DynMan}
on limit laws of degenerate symmetric statistics.
In this section we briefly review these results.

 Let $\clx$ be a Polish space and let  $\{X_n\}_{n=1}^{\infty}$ be a sequence of independent
identically distributed $\mathcal{X}$-valued random variables having common probability law $\nu$.
For $k=1,2,\dots$ let $L^2(\nu^{\otimes k})$ be the space of all real valued square
integrable functions on $(\mathcal{X}^k, \clb(\clx)^{\otimes k}, \nu^{\otimes k})$. Denote by $L^2_{sym}(\nu^{\otimes k})$ 
the subspace of symmetric functions, namely functions $\phi \in L^2(\nu^{\otimes k})$ such that for every permutation $\pi$ on $\{1, \cdots k\}$,
$$\phi(x_1, \cdots , x_k) = \phi(x_{\pi(1)}, \cdots , x_{\pi(k)}), \; \nu^{\otimes k} \mbox{ a.e }  (x_1,  \ldots x_k).$$
Given  $\phi_k \in L^2_{sym}(\nu^{\otimes k})$ define a symmetric statistic  $\sigma_k^n(\phi_k)$ as
\begin{align*}\sigma_k^n(\phi_k)&= \sum_{1\leq i_i<i_2\dots<i_k\le n}\phi_k(X_{i_1},\dots,X_{i_k})\,\,\,\mbox{for}\,\,n\geq k\\
&=0\;\;\;\;\;\;\;\;\;\;\;\;\;\;\;\;\;\;\;\;\;\;\;\;\;\;\;\;\;\;\;\;\;\;\;\;\;\;\;\;\;\;\;\;\;\;\;\mbox{for}\,\,n<k\,.
\end{align*}
In order to describe the asymptotic distributions of such statistics consider a Gaussian field
 $\{I_1(h);h\in L^2(\nu)\}$ such that 
$$\E(I_1(h)) = 0, \; \; \mathbb{E}(I_1(h)I_1(g)) = \langle h, g\rangle_{L^2(\nu)}\; \mbox{ for all } h, g \in L^2(\nu),$$ 
where $\langle \cdot, \cdot\rangle_{L^2(\nu)}$ denotes the inner product in $L^2(\nu)$.
For  $h \in L^2(\nu)$ define $\phi_k^h \in L^2_{sym}(\nu^{\otimes k})$, $k \ge 1$ as
\begin{align*}\phi_k^h(x_1, \ldots, x_k) = h(x_1) \cdots h(x_k),\; (x_1,  \ldots , x_k) \in \clx^k.
\end{align*}
We set
$\phi^h_0=1$.

The  multiple Wiener integral(MWI) of $\phi^h_k$, denoted as  $I_k(\phi^h_k)$, is defined through the following formula.
For $k\ge 1$
$$
I_k(\phi_k^h) = \sum_{j=0}^{\lfloor k/2 \rfloor} (-1)^j C_{k,j} \|h\|^{2j}_{L^2(\nu)} (I_1(h))^{k-2j}, \; \mbox{ where } C_{k,j} = \frac{k!}{(k-2j)!\, 2^j\, j!}, \; j= 0, \ldots , \lfloor k/2 \rfloor.$$
 The following representation gives an equivalent way to characterize MWI of $\phi_k^h$, $k \ge 1$.
\begin{align*}
	\sum_{k=0}^{\infty}\frac{t^k}{k!}I_k(\phi^h_k)=
	\exp\left (tI_1(h)-\frac{t^2}{2}\|h\|^2_{L^2(\nu)}\right),\; t \in \RR ,
\end{align*}
	where we set $I_0(\phi_0^h) = 1$.
We extend the definition of $I_k$ to the linear span of $\{\phi^h_k, h \in L^2(\nu)\}$ by linearity. It can be checked that for all $f$ in this linear span
\begin{equation}\E(I_k(f))^2 = k! \|f\|^2, \label{eq:isom}\end{equation}
	where on the right side $\| \cdot \|$ denotes the usual norm in $L^2(\nu^{\otimes k})$.
	Using this identity and standard denseness arguments, the definition of $I_k(f)$ can be extended to all $f \in L^2_{sym}(\nu^{\otimes k})$ and  the identity \eqref{eq:isom} holds for
	all $f \in L^2_{sym}(\nu^{\otimes k})$.
	The following theorem is taken from \cite{DynMan}.
\begin{theorem}[Dynkin-Mandelbaum \cite{DynMan}] 
	\label{thm:DM}
	Let $\{\phi_k\}_{k=1}^{\infty}$ be such that,
for each $k\geq 1$, $\phi_k \in L^2_{sym}(\nu^{\otimes k})$, and
\begin{align*}\int\phi_k(x_1,\dots,x_{k-1},x)\nu(dx)=0\,\,\, \mbox{for}\,\,\,\nu^{\otimes k-1} \mbox{ a.e.} \,(x_1,\dots,x_{k-1})\,.
\end{align*} Then
\begin{align*}
	\left( n^{-\frac{k}{2}}\sigma^n_k(\phi_k)\right)_{k\ge 1} \Rightarrow \left( \frac{1}{k!}I_k(\phi_k) \right)_{k\ge 1}
\end{align*} 
as a sequence of $\RR^{\infty}$ valued random variables.
\end{theorem}
\section{Girsanov Change of Measure}
\label{sec:gir}
For $N \in \N$, let $\bar \Om_N$, $P^N$, $V^i$, $i =0, \ldots , N$, $Y$, $\mu^N$ be as in Section \ref{sec:cltassu}. Also let $\mu = \Pi(Y_0, B^0, \bfn^0)$.
With these definitions \eqref{eq:mfield1}-\eqref{eq:mfield2}  are satisfied for $i=1, \ldots , N$; $\mu_s = \cll(X^i(s) \mid \clg^0) = \cll(X^i(s) \mid \clg^0_s)$,
$s \in [0,T]$, $i=1, \ldots , N$;
and 
$Y$ is $\{\clg^0_t\}$ adapted, where
$\clg^0_s = \sigma \{ Y_0,B^0_r, \bfn^0([0,r]\times A), r \le s, A \in \clb(\RR^m \times \RR_+ )\}$ and $\clg^0 = \clg^0_T$.

In addition to the above processes, define $Y^N$ as the unique solution of the following equation
\begin{equation}
	\label{eq:yneq}
	Y^N_t = Y_0 + \int_0^t b_0(Y^N_s, \mu^N_s) ds + \int_0^t \sigma_0(Y^N_s, \mu^N_s) dB^0_s + \int_{\XX^0_t} k 1_{[0, d_0(k)]}(u)  d\bfn^0,
\end{equation}
where $\mu^N_s = \frac{1}{N} \sum_{i=1}^N \delta_{X^i_s}$.

  Let for $i = 1, \ldots , N$, $u \in \RR_+$, $h \in \RR^d$ and $s \in [0,T]$
$$R^i = (X^i, Y), \; R^{N,i} = (X^i, Y^N), \; \beta^{N,i}_s = b(R^{N,i}_s,  \mu^N_s) - b(R^i_s,  \mu_s),$$
$$\bfd^{N,i}_s(h) = d(R^{N,i}_s, \mu^N_s, h), \; \bfd^{i}_s(h) = d(R^{i}_s, \mu_s, h), $$
$$e^{N,i}_s(h) = \bfd^{N,i}_s(h) - \bfd^{i}_s(h), \; r^{N,i}_s(u,h) = 1_{[0,\bfd^{i}_s(h)]}(u) \log \frac{\bfd^{N,i}_s(h)}{\bfd^{i}_s(h)} .$$
Write $B^i = (W^i, \tilde W^i)$, where $W^i, \tilde W^i$ are independent $d$ and $p$ dimensional Brownian motions respectively.
Define $\{H^N(t)\}$ as
 $$H^N(t)=\exp\left(J^{N,1}(t)+ J^{N,2}(t)\right)$$ where
\begin{align*}
J^{N,1}(t)=\sum_{i=1}^N\left(\int_0^t\beta^{N,i}_s\cdot dW^i_s-\frac{1}{2}\int_0^t\|\beta^{N,i}_s\|^2ds\right)
\end{align*}and
\begin{align*}J^{N,2}(t)=\sum_{i=1}^N\left(\int_{\mathbb{X}_t}r^{N,i}_{s-}(u,h)d\bfn^i-\int_{[0,t]\times \mathbb{R}^d}e^{N,i}_s(h)\gamma(dh)ds\right)\,.
\end{align*}
Letting for $t \in [0, T]$, $\bar \clf^N_t = \sigma\{V^i(s), 0 \le s \le t, \; i = 0, \ldots , N\}$, we see that
 $\{H^N_t\}$ is a $\bar \clf^N_t$ martingale under $P^N$. 
Define a new probability measure $Q^N$ on $\bar \Omega_{N}$ by 
$$\frac{dQ^N}{dP^N}=H^N(T)\,.$$ 
Expected values under $P^N$ and $Q^N$ will be denoted as $\EE_{P^N}$ and $\EE_{Q^N}$ respectively.

By Girsanov's theorem,  $\{(X^1,\ldots, X^N, Y^N,V^0)\}$  has the same probability law under $Q^N$ as 
 $\{(Z^{N,1},\ldots, \\Z^{N,N}, U^N, \bar V^0)\}$ (defined in \eqref{eq:nparteq} - \eqref{eq:comminp} and above Theorem \ref{thm:clt}) under $P^N$. 
Thus in order to prove the theorem it suffices to show
that 
\begin{align*} &\lim_{N\rightarrow \infty}\mathbb{E}_{Q^N}\exp\biggl(i\bigl\{\sqrt{N}\bigl(\frac{1}{N}\sum_{j=1}^N\varphi(X^j)-m_{\varphi}(V^0)\bigr)\bigr\}\biggr)\\
&\quad 	=\int_{\Omega_m}\exp\biggl(-\frac{1}{2}(\sigma_{\omega_0}^{\varphi})^2\biggr)P_0(d\omega_0)\,,
\end{align*}
which is equivalent to showing
\begin{align}&\lim_{N\rightarrow \infty}\mathbb{E}_{P^N}\exp\biggl(i\bigl\{\sqrt{N}\bigl(\frac{1}{N}\sum_{j=1}^N\varphi(X^j)-m_{\varphi}(V^0)\bigr)\bigr\}+J^{N,1}(T)+J^{N,2}(T)\biggr)\nonumber\\
	&\quad =\int_{\Omega_m}\exp\biggl(-\frac{1}{2}(\sigma_{\omega_0}^{\varphi})^2\biggr)P_0(d\omega_0)\,.\label{eq:eq209}
\end{align}
This will be shown in Section \ref{sec:compprf}.  We begin with some estimates.

\section{Estimating $Y^N-Y$.}
\label{sec:estyny}
The following lemma is immediate from the fact that, under $P^N$, $\{X^j\}_{j\in \NN}$ are iid, conditionally on $\clg^0$.  We omit the proof.
\begin{lemma}
	\label{lem:elem}
	For each $l \in \NN$, there exists $\vartheta_{l} \in (0, \infty)$ such that for all $t \in [0,T]$
	$$\sup_{\|f\|_{\infty} \le 1} \mathbb{E}_{P^N} | \langle f, \mu_t-\mu_t^N\rangle|^l \le \frac{\vartheta_l}{N^{l/2}}.$$
\end{lemma}
We now introduce a condition on the coefficients $b_0$ and $\sigma_0$.  Write $\sigma_0 = (\sigma_0^1, \cdots , \sigma_0^m)$, where each $\sigma_0^i$ is a function with values in $\RR^m$.
 Denote by $\hat \clj$ the collection of all real functions $f$ on  $\RR^{m+d}$ that are bounded by $1$
 and are such that $x \mapsto  f(\tilde y , x)$ is continuous for all 
 $\tilde y \in \RR^{m}$.  We say a function $\psi: \RR^m \times \clp(\RR^d) \to \RR^m$ is in class $\cls_1$ if there exist $\hat c_{\psi}\in (0, \infty)$, a finite subset
$\hat \clj_{\psi}$ of $\hat \clj$, continuous and bounded functions
$\psi^{(1)},\psi^{(2)}$ from $\RR^{m} \times \clp(\RR^d)$ to $\RR^{m\times m}$ and $\RR^{m} \times \clp(\RR^d)\times \RR^d$ to $\RR^m$ respectively;
 		and $\theta_{\psi}: \RR^{m}\times \RR^{m}\times \clp(\RR^d) \times \clp(\RR^d) \to \RR^m$
such that for all $y, y' \in \RR^{m}$ and $\nu, \nu' \in \clp(\RR^d)$ 
\begin{equation}\label{eq:smooth1}
\psi(y',\nu') - \psi(y,\nu) =  \psi^{(1)}(y,\nu)(y'-y)  + \langle  \psi^{(2)}(y,\nu,\cdot), (\nu'-\nu)  \rangle +
\theta_{\psi}(y,y',\nu,\nu'),\end{equation}
where
\begin{equation}\label{eq:eq850}
	\|\theta_{\psi}(y,y',\nu,\nu')\| \le c_{\psi} \left ( \|y'-y\|^2 + \max_{f \in \hat \clj_{\psi}} 
	|\langle  f(y, \cdot), (\nu'-\nu)\rangle|^2\right),
\end{equation}
Furthermore,
		\begin{equation}\label{eq:eq808}
		\|\psi(y,\nu) - \psi(y',\nu')\| \le   c_{\psi} \left ( \|y-y'\| + \max_{f \in \hat \clj_{\psi}} 
		|\langle  f(y, \cdot), (\nu-\nu')\rangle|\right).\end{equation}
\begin{condition}
	\label{cond:clt2}
	The functions $b_0$ and $\sigma_0^i$, $i=1, \cdots , m$ are in $\cls_1$.
\end{condition}
\begin{lemma}
	\label{lem:ynybd}
	Suppose that Conditions \ref{cond:main}, \ref{cond:clt1} and \ref{cond:clt2} hold.  Then for each $l \in \NN$, there exists
	a $\tilde \vartheta_l \in (0, \infty)$, such that for all $t \in [0,T]$
	$$\mathbb{E}_{P^N} \|Y^N_t-Y_t\|^l \le \frac{\tilde \vartheta_l}{N^{l/2}}.$$
\end{lemma}
{\bf Proof.}  
Fix $l \in \NN$ and $t\in [0,T]$.  By standard martingale inequalities and property \eqref{eq:eq808} for $\psi = b_0, \sigma_0^i$, $i=1, \cdots , m$, we have that for some $k_l \in (0,\infty)$
$$
\mathbb{E}_{P^N} \|Y^N_t-Y_t\|^l \le k_l \mathbb{E}_{P^N} \int_0^t \|Y^N_s-Y_s\|^l ds + k_l \mathbb{E}_{P^N} \int_0^t \max_{f \in \hat \clj_{1}} 
|\langle  f(Y_s, \cdot), (\mu_s^N-\mu_s)\rangle|^l ds,$$
where $\hat \clj_{1} = (\hat \clj_{b_0})\cup(\cup_{i=1}^l\hat \clj_{\sigma_0^i})$.  The result is now immediate from Gronwall's lemma and Lemma \ref{lem:elem}.
\qed

The following lemma follows on using classical existence/uniqueness results for SDE and an application of Ito's formula.
We will use the following notation
 \begin{equation}\label{eq:clvnot}
	\clv = (Y,B^0, \mu),\; \clz = (Y,\mu),\; \clz^N = (Y^N, \mu^N).
\end{equation}
\begin{lemma}
	\label{lem:repnyny}
	Suppose that Conditions \ref{cond:main}, \ref{cond:clt1} and \ref{cond:clt2} hold. For $t \in [0,T]$
	$$Y^N_t - Y_t = \frac{1}{N} \sum_{j=1}^N \bfs_{0,t}(X^j_{[0,t]}, \clv_{[0,t]}) + \clt^N(t),$$
	where
	\begin{align*}
		\clt^N(t) &= \cle_t \int_0^t \cle_s^{-1} \theta_{b_0}(\clz_s, \clz^N_s) ds + 
		\sum_{k=1}^m 
		\cle_t \int_0^t \cle_s^{-1} \theta_{\sigma_0^k}(\clz_s, \clz^N_s) dB^{0,k}_s\\
		&\quad - \sum_{k=1}^m 
		\cle_t \int_0^t \cle_s^{-1} \sigma_0^{k,(1)}(\clz_s) \theta_{\sigma_0^k}(\clz_s, \clz^N_s) ds,
	\end{align*}
	\begin{align*}
		\bfs_{0,t}(X^j_{[0,t]}, \clv_{[0,t]}) &= \cle_t \int_0^t \cle_s^{-1} b_0^{(2),c}(\clz_s, X^j_s) ds + \sum_{k=1}^m 
		\cle_t \int_0^t \cle_s^{-1} \sigma_0^{k,(2),c}(\clz_s, X^j_s) dB^{0,k}_s\\
		&\quad - \sum_{k=1}^m 
		\cle_t \int_0^t \cle_s^{-1} \sigma_0^{k,(1)}(\clz_s) \sigma_0^{k,(2),c}(\clz_s, X^j_s) ds,
	\end{align*}
	$\{\cle_t\}$ solves the $m\times m$ dimensional SDE
	$$
	\cle_t = I + \int_0^t b_0^{(1)}(\clz_s) \cle_s ds + \sum_{k=1}^m \int_0^t \sigma_0^{k,(1)}(\clz_s) \cle_s dW^k_s,$$
	$$	b_{0}^{(2),c}(y,\nu, \tilde x) = b_{0}^{(2)}(y,\nu, \tilde x) - \int_{\RR^d}  b_{0}^{(2)}(y,\nu,   x') \nu(dx'), \; (y,\nu , \tilde x) \in \RR^m \times \clp(\RR^d) \times \RR^d$$
	and $\sigma_0^{k,(2),c}$ is defined similarly.
\end{lemma}
{\bf Proof.} Using \eqref{eq:smooth1} with $\psi = b_0, \sigma_0^i$, $i=1, \cdots , m$, we have that 
\begin{align*}
	Y^N_t - Y_t &= \int_0^t (b_0(\clz^N_s)-b_0(\clz_s)) ds + \int_0^t (\sigma_0(\clz^N_s)-\sigma_0(\clz_s)) dB^0(s)\\
	&= \int_0^t \left(b_0^{(1)}(\clz_s) (Y^N_s - Y_s) + \langle b_0^{(2)}(\clz_s,\cdot), \mu^N_s-\mu_s\rangle + \theta_{b_0}(\clz_s, \clz^N_s) \right) ds\\
	&\quad + \sum_{k=1}^m \int_0^t \left(\sigma_0^{k,(1)}(\clz_s) (Y^N_s - Y_s) + \langle \sigma_0^{k,(2)}(\clz_s,\cdot), \mu^N_s-\mu_s\rangle + \theta_{\sigma_0^k}(\clz_s, \clz^N_s) \right) dB^{0,k}_s.\end{align*}
A standard application of Ito's formula shows
\begin{align*}
	Y^N_t - Y_t &= \cle_t \int_0^t \cle_s^{-1} \left(\langle b_0^{(2)}(\clz_s,\cdot), \mu^N_s-\mu_s\rangle + \theta_{b_0}(\clz_s, \clz^N_s)\right) ds\\
	&\quad + \sum_{k=1}^m \cle_t \int_0^t \cle_s^{-1} \left(\langle \sigma_0^{k,(2)}(\clz_s,\cdot), \mu^N_s-\mu_s\rangle + \theta_{\sigma_0^k}(\clz_s, \clz^N_s)\right) dB^{0,k}_s\\
	&\quad - \sum_{k=1}^m \cle_t \int_0^t \cle_s^{-1} \sigma_0^{k,(1)}(\clz_s) \left(\langle \sigma_0^{k,(2)}(\clz_s,\cdot), \mu^N_s-\mu_s\rangle + \theta_{\sigma_0^k}(\clz_s, \clz^N_s)\right) ds.
\end{align*} 
The result now follows on rearranging terms and noting that
$$\langle b_0^{(2)}(\clz_s,\cdot), \mu^N_s-\mu_s\rangle = \frac{1}{N} \sum_{j=1}^N b_0^{(2),c}(\clz_s, X^j_s),\; \langle \sigma_0^{k,(2)}(\clz_s,\cdot), \mu^N_s-\mu_s\rangle = \frac{1}{N} \sum_{j=1}^N\sigma_0^{k,(2),c}(\clz_s, X^j_s).$$
\qed

The following lemma follows on using the boundedness of coefficients, an application of Gronwall's lemma, Holder's inequality, Lemmas \ref{lem:elem} and \ref{lem:ynybd} and properties of $\theta_{\psi}$ 
for $\psi$ in class $\cls_1$.  The proof is omitted.
Let
$\bfh^j_t = \bfs_{0,t}(X^j_{[0,t]}, \clv_{[0,t]})$, $t \in [0,T]$.
\begin{lemma}
	\label{lem:etbds}
	Suppose that Conditions \ref{cond:main}, \ref{cond:clt1} and \ref{cond:clt2} hold.  Then for each $l \in \NN$
	$$
	\sup_{n\in \NN} \sup_{t\in [0,T]} \left ( \mathbb{E}_{P^N}\|\cle_t\|^l + \mathbb{E}_{P^N}\|\cle_t^{-1}\|^l\right) < \infty$$
	and there exists $\varpi \in (0, \infty)$ such that for all $t\in [0,T]$
	$$\mathbb{E}_{P^N}\|\frac{1}{N} \sum_{j=1}^N \bfh^j_t\|^2 \le \frac{\varpi}{N},\; \mathbb{E}_{P^N}\|\clt^N(t)\|^2 \le \frac{\varpi}{N^2}.$$
\end{lemma}

\section{Proof of Theorem \ref{thm:clt}.}
\label{sec:proofmainth}

\subsection{Asymptotics of $J^{N,1}$.}
\label{sec:7.1}
In Lemmas \ref{lem:term1} and \ref{lem:term2} below we study the asymptotics of the first and second sums
in $J^{N,1}$ respectively.  For this we introduce an additional condition on the coefficient $b$.  Denote by $\clj$  the collection of all real functions $f$ on $\RR^{d+m+d}$  that are bounded by $1$
 and are such that $x \mapsto f(\tilde x, \tilde y , x)$  is continuous for all 
 $(\tilde x, \tilde y) \in \RR^{d+m}$.
\begin{condition}
	\label{cond:clt3}
	There exist $c_b \in (0, \infty)$; a finite subset $\clj_F$ of $\clj$; continuous and bounded  functions $b_2, b_3$ from $\RR^{d+m} \times \clp(\RR^d)$ to $\RR^{d\times m}$ and $\RR^{d+m} \times \clp(\RR^d)\times \RR^d$ to $\RR^d$ respectively;
			and $\theta_b: \RR^{d+m}\times \RR^{d+m}\times \clp(\RR^d) \times \clp(\RR^d) \to \RR^d$ such that for all $z=(x,y), z' = (x,y') \in \RR^{d+m}$ and $\nu, \nu' \in \clp(\RR^d)$ 
				$$
				b(z',\nu') - b(z,\nu) =  b_2(z,\nu)(y'-y)  + \langle b_3(z,\nu,\cdot), (\nu'-\nu)\rangle +
				\theta_b(z,z',\nu,\nu')$$
				and 
				\begin{equation}\label{eq:eq802}
					\|\theta_b(z,z',\nu,\nu')\| \le c_b \left ( \|y'-y\|^2 + \max_{f \in \clj_F} 
					|\langle  f(z, \cdot), (\nu'-\nu)\rangle|^2\right).\end{equation}
					\end{condition}
\begin{lemma}
	\label{lem:term1}
	Suppose that Conditions \ref{cond:main}, \ref{cond:clt1}, \ref{cond:clt2}, \ref{cond:clt3} hold.  For $N \in \N$,
	$$\sum_{i=1}^N \int_0^T \beta^{N,i}_s dW^i_s = \frac{1}{N}\sum_{i \neq j} \int_0^T b_3^c(R^i_s, \mu_s, X^j_s) dW^i_s
	+ \frac{1}{N}\sum_{i \neq j} \int_0^T b_2(R^i_s, \mu_s) \bfh^j_s dW^i_s + \clr^N_1,$$
	where  $\clr^N_1$ 
	converges to $0$ in probability, where 
		$$b_{3}^c(x,y,\nu, \tilde x) = b_{3}(x,y,\nu, \tilde x) - \int_{\RR^d}  b_{3}(x,y,\nu,  x') \nu(dx'), \; (x,y,\nu , \tilde x) \in \RR^{d+m} \times \clp(\RR^d) \times \RR^d.$$.
\end{lemma}
{\bf Proof.}  
By Condition \ref{cond:clt3} it follows that, for $s \in [0,T]$,
\begin{equation}
	\label{eq:eq750}
\beta^{N,i}_s = b_2(R^i_s, \mu_s)(Y^N_s - Y_s) + \langle  b_3(R^i_s, \mu_s, \cdot), (\mu^N_s - \mu_s)\rangle + \zeta^{N,i}_s,
\end{equation}
where $\zeta^{N,i}_s = \theta_b(R^i_s, R^{N,i}_s, \mu_s, \mu^N_s)$.  Next, from \eqref{eq:eq802} we have
\begin{align*}
	\EE_{P^N} \left(\sum_{i=1}^N \int_0^T \zeta^{N,i}_s dW^i_s\right)^2 &= \sum_{i=1}^N \int_0^T \EE_{P^N}\|\zeta^{N,i}_s\|^2 ds \\
	&\le \kappa_1 \sum_{i=1}^N \int_0^T \EE_{P^N}(\|Y^N_s - Y_s\|^4)ds\\
	&+ \kappa_1 \sum_{i=1}^N \sum_{f \in \clj_F} \int_0^T \EE_{P^N}(\langle f(R^i_s, \cdot), (\mu^N_s - \mu_s)\rangle)^4ds.
\end{align*}
Since $f$ is bounded by $1$; $\clj_F$ is a finite collection; and conditionally on $\clg^0$, $X^i$ are i.i.d., the second term on the right side using Lemma \ref{lem:elem} can be bounded by 
$\kappa_2/N$ for some $\kappa_2 \in (0, \infty)$.
Also, from Lemma \ref{lem:ynybd} the first term converges to $0$.
	Combining the above observations we have, as $N \to \infty$,
	\begin{equation}
		\label{eq:eq815}
		\sum_{i=1}^N \int_0^T \zeta^{N,i}_s dW^i_s  \to  0, \mbox{ in probability }.
		\end{equation}
		Now consider the second term in \eqref{eq:eq750}:
		\begin{align}
		&\sum_{i=1}^N \int_0^T	\langle  b_3(R^i_s, \mu_s, \cdot), (\mu^N_s - \mu_s)\rangle dW^i_s \nonumber\\
		&\quad = \frac{1}{N} \sum_{i=1}^N \int_0^T b_3^c(R^i_s, \mu_s, X^i_s) dW^i_s +
		\frac{1}{N} \sum_{i\neq j} \int_0^T b_3^c(R^i_s, \mu_s, X^j_s) dW^i_s. \label{eq:eq856}
		\end{align}
		
Using the boundedness of $b_3$ it follows that,
\begin{equation}
	\label{eq:eq857}
\frac{1}{N} \sum_{i=1}^N \int_0^T b_3^c(R^i_s, \mu_s, X^i_s) dW^i_s \to  0, \mbox{ in probability }.\end{equation}
Finally consider the first term in \eqref{eq:eq750}.  
From Lemma \ref{lem:repnyny}, for $t \in [0,T]$,
\begin{align*}
	\sum_{i=1}^N \int_0^t b_2(R^i_s, \mu_s)(Y^N_s - Y_s) dW^i_s &= 
	\frac{1}{N} \sum_{i=1}^N \int_0^t b_2(R^i_s, \mu_s) \bfh^i_s dW^i_s + \frac{1}{N} \sum_{i\neq j} \int_0^t b_2(R^i_s, \mu_s) \bfh^j_s dW^i_s\\
	&+ \sum_{i=1}^N \int_0^t b_2(R^i_s, \mu_s) \clt^N(s) dW^i_s.
\end{align*}
The first term on the right side converges  to $0$ in probability since $b_2, b_{0}^{(i)}, \sigma_0^{k, (i)}$  are bounded.
Also, using the boundedness of $b_2$ and Lemma \ref{lem:etbds}, the third term converges to $0$ in probability.
Result now follows on combining the above observation with \eqref{eq:eq815}, \eqref{eq:eq856} and \eqref{eq:eq857}. \qed

For the next lemma we will need some notation.
Define functions $\bfs_{1,t}, \bfs^c_{1,t}$ from $\RR^d \times \DD_{\RR^{2d+m} \times \clp(\RR^d)}[0,t]$ to $\RR$ as follows:
For $(x,  x^{(1)}_{[0,t]}, x^{(2)}_{[0,t]}, y_{[0,t]}, w_{[0,t]}, \nu_{[0,t]}) \equiv (x, \zeta_{[0,t]}) \in \RR^d \times \DD_{\RR^{2d+2m} \times \clp(\RR^d)}[0,t]$
\begin{align*}
	\bfs_{1,t}(x,  \zeta_{[0,t]}) &=
	b_2(x, y_t, \nu_t)\bfs_{0,t}(\zeta^{(1)}_{[0,t]})\cdot b_2(x, y_t, \nu_t)\bfs_{0,t}(\zeta^{(2)}_{[0,t]})  \\
	\bfs_{1,t}^c(x, \zeta_{[0,t]}) &= \bfs_{1,t}(x, \zeta_{[0,t]})
	- \bfm_{1,t}(\zeta_{[0,t]})
	\end{align*}
where $\zeta^{(i)}_{[0,t]} = (x^{(i)}_{[0,t]}, y_{[0,t]}, w_{[0,t]}, \nu_{[0,t]})$,
and the function $\bfm_{1,t}$ from $\DD_{\RR^{2d+2m} \times \clp(\RR^d)}[0,t]$ to $\RR$ is defined as
$$\bfm_{1,t}(\zeta_{[0,t]}) = \int_{\RR^d} \bfs_{1,t}(x', \zeta_{[0,t]}) \nu_t(dx').$$

Next, define for $t \in [0,T]$, functions $\bfs_{2,t}, \bfs^c_{2,t}$ from $\RR^d \times \DD_{\RR^{2d+2m} \times \clp(\RR^d)}[0,t]$ to $\RR$ as follows:
\begin{align*}
	\bfs_{2,t}(x, \zeta_{[0,t]}) &=
 2b_2(x, y_t, \nu_t) \bfs_{0,t}(\zeta^{(1)}_{[0,t]}) \cdot b_3^c(x, y_t, \nu_t,  x^{(2)}_t) \\
	\bfs_{2,t}^c(x, \zeta_{[0,t]}) &= \bfs_{2,t}(x, \zeta_{[0,t]})
	- \bfm_{2,t}(\zeta_{[0,t]})
\end{align*}
where the function $\bfm_{2,t}$ from $\DD_{\RR^{2d+2m} \times \clp(\RR^d)}[0,t]$ to $\RR$ is defined as
$$\bfm_{2,t}(\zeta_{[0,t]}) = 
\frac{1}{2} \sum_{i,j \in \{1,2\}, i \neq j}
\int_{\RR^d} \bfs_{2,t}(x', x^{(i)}_{[0,t]}, x^{(j)}_{[0,t]}, y_{[0,t]}, w_{[0,t]}, \nu_{[0,t]}) \nu_t(dx').$$

Also, define functions $\bfs_3, \bfs_3^c$ from $\RR^{3d+m} \times \clp(\RR^d)$ to $\RR$ as follows: For $(x,x^{(1)}, x^{(2)}, y,\nu) \in \RR^{3d+m} \times \clp(\RR^d)$
\begin{align*}
\bfs_3(x,x^{(1)}, x^{(2)},y,\nu) &= b_3^c(x,y,\nu ,x^{(1)})\cdot b_3^c(x,y,\nu , x^{(2)}), \\
\bfs_3^c(x,x^{(1)}, x^{(2)},y,\nu) &= \bfs_3(x, x^{(1)}, x^{(2)},y,\nu) - \bfm_3(x^{(1)}, x^{(2)}, y,\nu),
\end{align*}
where $\bfm_3$  from $\RR^{2d+m} \times \clp(\RR^d)$ to $\RR$ is defined as
$$\bfm_3(x^{(1)}, x^{(2)}, y,\nu) = \int \bfs_3(x',x^{(1)}, x^{(2)},y,\nu) \nu(dx').$$

Finally, define $\bfm_t$ from $\DD_{\RR^{2d+2m} \times \clp(\RR^d)}[0,t]$ to $\RR$ as follows.
$$
\bfm_{t}(\zeta_{[0,t]}) = \sum_{i=1}^2\bfm_{i,t}(\zeta_{[0,t]})
+ \bfm_3(x^{(1)}_t, x^{(2)}_t, y_t,\nu_t).$$
Recall the process $\clv$ from \eqref{eq:clvnot}.
\begin{lemma}
	\label{lem:term2}
	For $N \in \N$,
	\begin{align}
		\sum_{i=1}^N \int_0^T \|\beta^{N,i}_s\|^2 ds
		&= \frac{1}{N} \sum_{j \neq k} \int_0^T \bfm_{t}(X^j_{[0,t]}, X^k_{[0,t]}, \clv_{[0,t]}) dt\nonumber\\
		&\quad + \frac{1}{N} \sum_{j=1}^N \int_0^T \bfm_{t}(X^j_{[0,t]}, X^j_{[0,t]}, \clv_{[0,t]}) dt
			+   \clr^N_2, \label{eq:eqquadform}
\end{align}
	where  $\clr^N_2$ 
	converges to $0$ in probability.
\end{lemma}
{\bf Proof.}
For $N \in \N$, $i = 1, \ldots N$ and $s \in [0,T]$
\begin{align}
	\|\beta^{N,i}_s\|^2 &= \|b(R^{N,i}_s,  \mu^N_s) - b(R^i_s,  \mu_s)\|^2\nonumber\\
	&= \|b_2(R^i_s,\mu_s)(Y^N_s-Y_s)  + \langle b_3(R^i_s,\mu_s,\cdot), (\mu^N_s-\mu_s)  \rangle +
	\theta_b(R^i_s,R^{N,i}_s,\mu_s,\mu^N_s)\|^2\nonumber\\
	&= \|b_2(R^i_s,\mu_s)(Y^N_s-Y_s)\|^2 + \|\langle  b_3(R^i_s,\mu_s,\cdot), (\mu^N_s-\mu_s)  \rangle\|^2\nonumber\\
	&\quad + \|\theta_b(R^i_s,R^{N,i}_s,\mu_s,\mu^N_s)\|^2 +
	2 b_2(R^i_s,\mu_s)(Y^N_s-Y_s) \cdot \langle  b_3(R^i_s,\mu_s,\cdot), (\mu^N_s-\mu_s) \rangle + \clt^{N,i}_1(s),\nonumber\\
	\label{eq:eq323d}
\end{align}
where the term $\clt^{N,i}_1(s)$ consists of the remaining two crossproduct terms.
Using \eqref{eq:eq802} and Lemmas \ref{lem:elem} and \ref{lem:ynybd},   we see that
\begin{equation}
	\label{eq:eq323bnew}
	\sum_{i=1}^N \int_0^T  \|\theta_b(R^i_s,R^{N,i}_s,\mu_s,\mu^N_s)\|^2ds \to 0 \mbox{ in probability as } N \to \infty .
\end{equation}
Similar estimates show that
\begin{equation}
	\label{eq:eq323c}
	\sum_{i=1}^N \int_0^T |\clt^{N,i}_1(s)|ds \to 0 \mbox{ in probability as } N \to \infty .
\end{equation}
Next, using Lemma \ref{lem:repnyny}, we have
\begin{align*}
\|b_2(R^i_s,\mu_s)(Y^N_s-Y_s)\|^2  &= \| \frac{1}{N} \sum_{j=1}^N b_2(R^i_s,\mu_s)\bfh^j_s \|^2 \\
&\quad + 
	\|b_2(R^i_s,\mu_s)\clt^N_s\|^2  +  \clt^{N,i}_2(s),
\end{align*}
where $\clt^{N,i}_2(s)$ is the corresponding crossproduct term.
Making use of Lemma \ref{lem:etbds} we can bound $\EE_{P^N,\clg^0}|\clt^{N,i}_2(s)|$ by
$\frac{\kappa_1}{N^{3/2}}$ for some $\kappa_1 > 0$ that does not depend on $s \in [0,T]$ and $i,N \in \N$.
Similarly, the  expected value of the second term in the above display can be bounded by $\frac{\kappa_2}{N^2}$ for some $\kappa_2 > 0$.  Thus 
\begin{equation}
	\label{eq:eq320}
	\sum_{i=1}^N \int_0^T \|b_2(R^i_s,\mu_s)(Y^N_s-Y_s)\|^2 ds = \sum_{i=1}^N \int_0^T \| \frac{1}{N} \sum_{j=1}^N b_2(R^i_s,\mu_s)\bfh^j_s \|^2 ds + \tilde \clr^N_1,
\end{equation}
where $\tilde \clr^N_1\to 0$ in probability as $N \to \infty$.

Recalling  the definition of $s_{1,t}$
\begin{align*}
	\int_0^T \sum_{i=1}^N \| \frac{1}{N} \sum_{j=1}^N b_2(R^i_s,\mu_s)\bfh^j_s \|^2  ds 
	& = \frac{1}{N^2} \sum_{i,j,k} \int_0^T b_2(R^i_s,\mu_s)\bfh^j_s \cdot b_2(R^i_s,\mu_s)\bfh^k_s ds \\
	&= \frac{1}{N^2} \sum_{i,j,k} \int_0^T \bfs_{1,t}(X^i_t,X^j_{[0,t]},X^k_{[0,t]}, \clv_{[0,t]}) dt.
\end{align*}
The above expression can be written as
\begin{align}
	 &\frac{1}{N^2} \sum_{i,j,k} \int_0^T \bfs^c_{1,t}(X^i_t,X^j_{[0,t]},X^k_{[0,t]}, \clv_{[0,t]}) dt
	+ \frac{1}{N} \sum_{j\neq k} \int_0^T \bfm_{1,t}(X^j_{[0,t]},X^k_{[0,t]}, \clv_{[0,t]}) dt\nonumber\\
	&\quad + \frac{1}{N} \sum_{j} \int_0^T \bfm_{1,t}(X^j_{[0,t]},X^j_{[0,t]}, \clv_{[0,t]}) dt .
\label{eq:eq1153}
	\end{align}	
	 From the boundedness of $\bfs^c_{1,t}$, conditional independence of $X^i, X^j, X^k$ for distinct indices $i,j,k$ and the fact
	that for all $(x, x^{(1)}_{[0,t]}, x^{(2)}_{[0,t]}, v_{[0,t]}) \in \RR^d \times \DD_{\RR^{2d+2m} \times \clp(\RR^d)}[0,t]$
	\begin{align*}
	\EE_{P^N} \bfs_{1,t}^c(X^i_t,x^{(1)}_{[0,t]},x^{(2)}_{[0,t]}, v_{[0,t]})	
	&= \EE_{P^N} \bfs_{1,t}^c(x,X^i_{[0,t]},x^{(2)}_{[0,t]}, v_{[0,t]})\\
	&= \EE_{P^N} \bfs_{1,t}^c(x,x^{(1)}_{[0,t]},X^j_{[0,t]}, v_{[0,t]}) = 0,
		\end{align*}
		it follows  that the first term in \eqref{eq:eq1153} 
		converges to $0$ in probability.

Next, 
$$\|\langle  b_3(R^i_s,\mu_s,\cdot), (\mu^N_s-\mu_s) \rangle\|^2
= \frac{1}{N^2} \sum_{j,k} b_3^c(R^i_s,\mu_s, X^{j}_s)\cdot b_3^c(R^i_s,\mu_s, X^{k}_s).$$
Thus
\begin{align*}
\sum_{i=1}^N \int_0^T \|\langle  b_3(R^i_s,\mu_s,\cdot), (\mu^N_s-\mu_s) \rangle\|^2 ds 
&=
\frac{1}{N^2} \sum_{i, j,k} \int_0^T b_3^c(R^i_s,\mu_s, X^{j}_s)\cdot b_3^c(R^i_s,\mu_s, X^{k}_s) ds  \\
&= \frac{1}{N^2} \sum_{i, j,k} \int_0^T \bfs_3(X^i_s,X^{j}_s, X^{k}_s, Y_s, \mu_s) ds  \\
\end{align*}
The above expression can be rewritten as
\begin{align}
 &\frac{1}{N^2} \sum_{i, j,k} \int_0^T \bfs_3^c(X^i_s,X^{j}_s, X^{k}_s, Y_s, \mu_s) ds 
 +  \frac{1}{N} \sum_{j\neq k} \int_0^T \bfm_3(X^{j}_s, X^{k}_s, Y_s, \mu_s) ds \nonumber \\
&\quad +  \frac{1}{N} \sum_{j} \int_0^T \bfm_3(X^{j}_s, X^{j}_s, Y_s, \mu_s) ds .
 \label{eq:eq947}
\end{align}
As before, the first term in \eqref{eq:eq947} converges  to $0$ in probability.

Finally we consider the crossproduct term in \eqref{eq:eq323d}:
\begin{align*} 
&\sum_{i=1}^N b_2(R^i_s,\mu_s)(Y^N_s-Y_s) \cdot \langle  b_3(R^i_s,\mu_s,\cdot), (\mu^N_s-\mu_s)  \rangle\\
&\quad = \frac{1}{N^2} \sum_{i,j,k} b_2(R^i_s, \mu_s)\bfh^j_s \cdot b_3^c(R^i_s, \mu_s, X^k_s) \\
&\quad + \frac{1}{N} \sum_{i,k} b_2(R^i_s, \mu_s)\clt^N(s) \cdot b_3^c(R^i_s, \mu_s, X^k_s)\\
&\quad \equiv \clt_3^N(s) + \clt_4^N(s)
\end{align*}
where the equality follows from Lemma \ref{lem:repnyny}.
Using Lemma \ref{lem:etbds}
we see that
$\int_0^T \clt^N_4(s) ds$  converges to $0$ in probability as $N \to \infty$.
For the  term $\clt_3^N(s)$
\begin{align*} 
\frac{2}{N^2} \sum_{i,j,k} \int_0^T b_2(R^i_s, \mu_s)\bfh^j_s \cdot b_3^c(R^i_s, \mu_s, X^k_s) ds
&= \frac{1}{N^2} \sum_{i,j,k} \int_0^T \bfs^c_{2,t}(X^i_t,X^j_{[0,t]},X^k_{[0,t]}, \clv_{[0,t]}) dt\\
&\quad + \frac{1}{N} \sum_{j\neq k} \int_0^T \bfm_{2,t}(X^j_{[0,t]},X^k_{[0,t]}, \clv_{[0,t]}) dt\nonumber\\
&\quad + \frac{1}{N} \sum_{j} \int_0^T \bfm_{2,t}(X^j_{[0,t]},X^j_{[0,t]}, \clv_{[0,t]}) dt.
\end{align*}
The first term on the right side once more converges to $0$ in probability.
The result now follows on combining the above display with  \eqref{eq:eq323d}, \eqref{eq:eq323bnew}, \eqref{eq:eq323c}, \eqref{eq:eq320}, \eqref{eq:eq1153} and \eqref{eq:eq947}.
\qed

\subsection{Asymptotics of $J^{N,2}$.}
\label{sec:secjn2}
We now consider the term $J^{N,2}$.  Recall the constants $\epsilon, K$ from Condition \ref{cond:main}.

From Taylor's expansion, there exists  a $\kappa_1 \in (0, \infty)$ such that
for all $\alpha, \beta \in (\epsilon, K)$
$$\log\frac{\alpha}{\beta} = (\frac{\alpha}{\beta}-1) - \frac{1}{2} (\frac{\alpha}{\beta}-1)^2 + \vartheta(\alpha, \beta) (\frac{\alpha}{\beta}-1)^3,$$
where $|\vartheta(\alpha, \beta)| \le \kappa_1$.  Letting $\vartheta^{N,i}_s(h) = \vartheta(\bfd^{N,i}_s(h), \bfd^{i}_s(h))$, we get
$$
\log\frac{\bfd^{N,i}_s(h)}{\bfd^{i}_s(h)} = \left(\frac{\bfd^{N,i}_s(h)}{\bfd^{i}_s(h)}-1\right) - \frac{1}{2}\left(\frac{\bfd^{N,i}_s(h)}{\bfd^{i}_s(h)}-1\right)^2
+ \vartheta^{N,i}_s(h)\left(\frac{\bfd^{N,i}_s(h)}{\bfd^{i}_s(h)}-1\right)^3.$$
Thus
\begin{align}
\int_{\mathbb{X}_T}r^{N,i}_{s-}(u,h)d\bfn^i-\int_{[0,T]\times \mathbb{R}^d}e^{N,i}_s(h)\gamma(dh)ds
&= \int_{\mathbb{X}_T} 1_{[0,\bfd^{i}_{s-}(h)]}(u) \left(\frac{\bfd^{N,i}_{s-}(h)}{\bfd^{i}_{s-}(h)}-1\right) d \tilde\bfn^i \nonumber\\
&\quad -\frac{1}{2} \int_{\mathbb{X}_T} 1_{[0,\bfd^{i}_{s-}(h)]}(u) \left(\frac{\bfd^{N,i}_{s-}(h)}{\bfd^{i}_{s-}(h)}-1\right)^2 d \bfn^i \nonumber \\
&\quad + \int_{\mathbb{X}_T} 1_{[0,\bfd^{i}_{s-}(h)]}(u) \vartheta^{N,i}_{s-}(h)\left(\frac{\bfd^{N,i}_{s-}(h)}{\bfd^{i}_{s-}(h)}-1\right)^3 d \bfn^i,
\label{eq:eq547}
\end{align}
where $\tilde\bfn^i$ is the compensated PRM: $\tilde \bfn^i = \bfn^i - \bfnu$.
In the lemmas below we consider the three terms on the right side of \eqref{eq:eq547} separately.
We introduce the following condition on the coefficient $d$.
 Denote by $\tilde \clj$ the collection of all real functions $f$ on  $\RR^{d+m+2d}$ that are bounded by $1$
 and are such that  $x \mapsto f(\tilde x, \tilde y , \tilde h, x)$ is continuous for all 
 $(\tilde x, \tilde y, \tilde h) \in \RR^{d+m+d}$.
\begin{condition}
	\label{cond:clt4}
	There exist $c_d \in (0, \infty)$; a finite subset $\tilde \clj_F$ of $\tilde \clj$; continuous and bounded real functions $d_2, d_3$ from $\RR^{d+m+d} \times \clp(\RR^d)$ to $\RR^m$ and $\RR^{d+m+d} \times \clp(\RR^d)\times \RR^d$ to $\RR$ respectively;
			and $\theta_d: \RR^{d+m}\times \RR^{d+m}\times \clp(\RR^d) \times \clp(\RR^d) \times \RR^d\to \RR$ such that for all $z=(x,y), z' = (x,y') \in \RR^{d+m}$, $h \in \RR^d$ and $\nu, \nu' \in \clp(\RR^d)$ 
			$$
			d(z',\nu',h) - d(z,\nu,h) = (y'-y) \cdot d_2(z,h,\nu) + \langle d_3(z,h,\nu,\cdot), (\nu'-\nu) \rangle +
			\theta_d(z,z',\nu,\nu',h)$$
			and 
			\begin{equation}\label{eq:eq206}
				|\theta_d(z,z',\nu,\nu')| \le c_d \left ( \|y'-y\|^2 + \max_{f \in \tilde \clj_F} |\langle f(z, h, \cdot), (\nu'-\nu) \rangle|^2\right).\end{equation}
\end{condition}
Next let  $d_3^c$ from $\RR^{d+m} \times \RR^d\times  \clp(\RR^d) \times \RR^d$ to $\RR$ as
$$d_3^c(x,y,h, \nu, \tilde x) = d_3(x,y,h, \nu, \tilde x) - \int_{\RR^d} d_3(x,y,h, \nu,  x') \nu(dx').$$
\begin{lemma}\label{lem:lem550}
	For $N \in \N$
	\begin{align*}
		\sum_{i=1}^N \int_{\mathbb{X}_T} 1_{[0,\bfd^{i}_{s-}(h)]}(u) \left(\frac{\bfd^{N,i}_{s-}(h)}{\bfd^{i}_{s-}(h)}-1\right) d \tilde\bfn^i
		&= \frac{1}{N} \sum_{i\neq j} \int_{\mathbb{X}_T} 1_{[0,\bfd^{i}_{s-}(h)]}(u) \frac{1}{\bfd^{i}_{s-}(h)} d_3^c(R^i_{s-}, h, \mu_{s-}, X^j_{s-}) d \tilde\bfn^i\\
		&\quad + \frac{1}{N} \sum_{i\neq j} \int_{\mathbb{X}_T} 1_{[0,\bfd^{i}_{s-}(h)]}(u) \frac{\bfh^j_{s}}{\bfd^{i}_{s-}(h)} d_2(R^i_{s-}, h, \mu_{s-}) d \tilde\bfn^i\\
		&\quad + \clr^N_3,
\end{align*}
where  $\clr^N_3$ 
converges to $0$ in probability.
\end{lemma}
{\bf Proof.}
From Condition \ref{cond:clt4}
\begin{equation}
	\label{eq:eq628}
	\bfd^{N,i}_s(h) - \bfd^{i}_s(h) = (Y^N_s-Y_s)\cdot d_2(R^i_s, h, \mu_s) + 
	\langle  d_3(R^i_s, h, \mu_s, \cdot), (\mu^N_s-\mu_s) \rangle
	+ \theta_d(R^i_s,  R^{N,i}_s, \mu_s, \mu^N_s, h).
\end{equation}
Since 
\begin{equation}\label{eq:eq747}
\sum_{i=1}^N \theta_d^2(R^i_s,  R^{N,i}_s, \mu_s, \mu^N_s, h) \le 
2c_d^2 \sum_{i=1}^N \left ( \|Y_s - Y^N_s\|^4 + \max_{f \in \tilde \clj_F} 
| \langle  f(R^i_s, h, \cdot), (\mu^N_s-\mu_s)\rangle |^4\right),
\end{equation}
we have from \eqref{eq:eq206}, Lemma \ref{lem:ynybd} and Lemma \ref{lem:elem} that, as $N\to \infty$
\begin{equation}
	\label{eq:eq638}
	\sum_{i=1}^N \int_{\mathbb{X}_T} 1_{[0,\bfd^{i}_{s-}(h)]}(u)
	\frac{1}{\bfd^{i}_{s-}(h)} \theta_d(R^i_{s-},  R^{N,i}_{s-}, \mu_{s-}, \mu^N_{s-}, h)\; d \tilde\bfn^i  \to 0
	\mbox{ in probability }.
\end{equation}
Next consider the second term on the right side of \eqref{eq:eq628}.
\begin{align}
	&\sum_{i=1}^N 1_{[0,\bfd^{i}_{s-}(h)]}(u) \frac{1}{\bfd^{i}_{s-}(h)}
	\langle  d_3(R^i_{s-}, h, \mu_{s-}, \cdot), (\mu^N_{s-}-\mu_{s-}) \rangle d \tilde\bfn^i \nonumber \\
	&\quad =
	\frac{1}{N} \sum_{i=1}^N 1_{[0,\bfd^{i}_{s-}(h)]}(u) \frac{1}{\bfd^{i}_{s-}(h)} d_3^c(R^i_{s-}, h, \mu_{s-}, X^i_{s-})  d \tilde\bfn^i\nonumber\\
	&\quad + \frac{1}{N} \sum_{i\neq j} 1_{[0,\bfd^{i}_{s-}(h)]}(u) \frac{1}{\bfd^{i}_{s-}(h)} d_3^c(R^i_{s-}, h, \mu_{s-}, X^j_{s-})  d \tilde\bfn^i.
	\label{eq:eq648}
\end{align}
Since $\{\bfn^i\}_{i=1}^N$ are independent, as $N \to \infty$,
\begin{equation}
	\label{eq:eq649}
\frac{1}{N}\sum_{i=1}^N \int_{\mathbb{X}_T} 1_{[0,\bfd^{i}_{s-}(h)]}(u) \frac{1}{\bfd^{i}_{s-}(h)} d_3^c(R^i_{s-}, h, \mu_{s-}, X^i_{s-})\;
d \tilde\bfn^i \to 0 \mbox{ in probability}.	
\end{equation}
Finally consider the first term on the right side of \eqref{eq:eq628}.  Using Lemma \ref{lem:repnyny}
\begin{align}
\label{eq:eq657}
(Y^N_s - Y_s)\cdot d_2(R^i_s, h, \mu_s) &= \frac{1}{N} \sum_{j=1}^N \bfh^j_s \cdot d_2(R^i_s, h, \mu_s) \nonumber \\
&\quad + 
\clt^N_s \cdot d_2(R^i_s, h, \mu_s) ds\nonumber\\
& \equiv \frac{1}{N} \sum_{j=1}^N \bfh^j_s \cdot d_2(R^i_s, h, \mu_s) + \hat{\clt}^{N,i}(s).
\end{align}
Using Lemma \ref{lem:etbds} we see that, as $N\to \infty$,
\begin{equation}
	\label{eq:eq704}
\sum_{i=1}^N \int_{\mathbb{X}_T} 1_{[0,\bfd^{i}_{s-}(h)]}(u) \frac{1}{\bfd^{i}_{s-}(h)} \hat{\clt}^{N,i}(s-)\;
d \tilde\bfn^i \to 0 \mbox{ in probability. }
\end{equation}
For the first term on the right side of \eqref{eq:eq657} note that
\begin{align}
&\frac{1}{N}\sum_{i,j=1}^N	
\int_{\mathbb{X}_T} 1_{[0,\bfd^{i}_{s-}(h)]}(u) \frac{\bfh^j_s}{\bfd^{i}_{s-}(h)}\cdot d_2(R^i_{s-}, h, \mu_{s-}) \;
d \tilde\bfn^i \nonumber \\
&= \frac{1}{N}\sum_{i=1}^N	\int_{\mathbb{X}_T} 1_{[0,\bfd^{i}_{s-}(h)]}(u) \frac{\bfh^j_s}{\bfd^{i}_{s-}(h)}\cdot 
d_2(R^i_{s-}, h, \mu_{s-}) \;
d \tilde\bfn^i\nonumber\\
&\quad + \frac{1}{N}\sum_{i\neq j}	\int_{\mathbb{X}_T} 1_{[0,\bfd^{i}_{s-}(h)]}(u) \frac{\bfh^j_s}{\bfd^{i}_{s-}(h)}
\cdot d_2(R^i_{s-}, h, \mu_{s-}) \;
d \tilde\bfn^i.
\label{eq:eq711}
\end{align}
As before, using the independence of $\{\bfn^i\}_{i=1}^N$, as $N\to \infty$,
$$
\frac{1}{N}\sum_{i=1}^N \int_{\mathbb{X}_T} 1_{[0,\bfd^{i}_{s-}(h)]}(u) \frac{\bfh^j_s}{\bfd^{i}_{s-}(h)}\cdot d_2(R^i_{s-}, h, \mu_{s-}) \;
d \tilde\bfn^i \to 0 \mbox{ in probability. }
$$
The result follows on combining the above display with  \eqref{eq:eq628}, \eqref{eq:eq638}, \eqref{eq:eq648}, \eqref{eq:eq649}, \eqref{eq:eq657}, \eqref{eq:eq704}
and \eqref{eq:eq711}. \qed

We now consider the second term on the right side of \eqref{eq:eq547}.  As for the proof of Lemma \ref{lem:term2}, we will need some additional notation.
Define for $t \in [0,T]$, function $\bar\bfs_{1,t}$ from $\RR^d \times \DD_{\RR^{2d+m} \times \clp(\RR^d)}[0,t]\times \RR^d$ to $\RR$ as follows:
For $(x, x^{(1)}_{[0,t]}, x^{(2)}_{[0,t]}, y_{[0,t]}, w_{[0,t]}, \nu_{[0,t]},h) \equiv (x,\zeta_{[0,t]})\in \RR^d \times \DD_{\RR^{2d+2m} \times \clp(\RR^d)}[0,t]\times \RR^d$
\begin{align*}
	\bar\bfs_{1,t}(x, \zeta_{[0,t]},h) &=
	\frac{1}{d(x,y_{t-},\nu_{t-},h)}
	\prod_{i=1}^2 \bfs_{0,t}(\zeta^{(i)}_{[0,t]}) \cdot d_2(x, y_{t-}, h, \nu_{t-}),
	%
	%
	%
	\end{align*}
	where $\zeta^{(i)}_{[0,t]} = (x^{(i)}_{[0,t]}, y_{[0,t]}, w_{[0,t]}, \nu_{[0,t]})$.
	Also define the function $\bar\bfm_{1,t}$ from $\DD_{\RR^{2d+2m} \times \clp(\RR^d)}[0,t]\times \RR^d$ to $\RR$ as
	$$\bar\bfm_{1,t}(\zeta_{[0,t]},h) = 
	\int_{\RR^d} \bar\bfs_{1,t}(x', \zeta_{[0,t]},h) \nu_t(dx').$$
	Next, define for $t \in [0,T]$, function $\bar\bfs_{2,t}$ from $\RR^d \times \DD_{\RR^{2d+2m} \times \clp(\RR^d)}[0,t]\times \RR^d$ to $\RR$ as follows:
	\begin{align*}
		\bar \bfs_{2,t}(x, \zeta_{[0,t]},h) &=
	 	\frac{2 d_3^c(x,y_{t-},h,\nu_{t-}, x^{(2)}_{t-})}{d(x,y_{t-},\nu_{t-},h)}
		 \bfs_{0,t}(\zeta^{(1)}_{[0,t]})
		\cdot d_2(x,y_{t-}, h,\nu_{t-}).
	\end{align*}
	Also define the function $\bar\bfm_{2,t}$ from $\DD_{\RR^{2d+m} \times \clp(\RR^d)}[0,t]\times \RR^d$ to $\RR$ as
	$$\bar\bfm_{2,t}(\zeta_{[0,t]},h) = 
	\frac{1}{2} \sum_{i,j \in \{1,2\}, i \neq j}\int_{\RR^d} \bar\bfs_{2,t}(x', x^{(i)}_{[0,t]},x^{(j)}_{[0,t]}, y_{[0,t]}, w_{[0,t]},\nu_{[0,t]},h) \nu_t(dx').$$
	Define function $\bar\bfs_3$ from $\RR^{3d+m} \times \clp(\RR^d)\times \RR^d$ to $\RR$ as follows: For $(x,x^{(1)},x^{(2)},y,\nu, h) \in \RR^{3d+m} \times \clp(\RR^d)\times \RR^d $
	$$
	\bar\bfs_3(x,x^{(1)},x^{(2)}, y,\nu,h) = \frac{\prod_{i=1}^2 d_3^c(x,y,h,\nu,x^{(i)})}{d(x,y,\nu,h)}$$
	and let $\bar \bfm_3$ be the function from $\RR^{2d+m} \times \clp(\RR^d)\times \RR^d$ to $\RR$ defined as
	$$\bar \bfm_3(x^{(1)},x^{(2)}, y,\nu,h) = \int \bar\bfs_3(x',x^{(1)},x^{(2)},y,\nu,h) \nu(dx').$$
	Finally, define $\bar\bfm_t$ from $\DD_{\RR^{2d+m} \times \clp(\RR^d)}[0,t]\times \RR^d$ to $\RR$ as follows.
	$$
	\bar\bfm_{t}(\zeta_{[0,t]},h) = \sum_{i=1}^2\bar\bfm_{i,t}(\zeta_{[0,t]},h)
	+ \bar\bfm_3(x^{(1)}_t, x^{(2)}_t, y_t,\nu_t,h).$$
	Recall the process $\clv$ introduced in \eqref{eq:clvnot}.
\begin{lemma}
	\label{lem:lem910}
	For $N \in \NN$
	\begin{align}
	&\sum_{i=1}^N \int_{\mathbb{X}_T} 1_{[0,\bfd^{i}_{s-}(h)]}(u) \left(\frac{\bfd^{N,i}_{s-}(h)}{\bfd^{i}_{s-}(h)}-1\right)^2 d \bfn^i\nonumber\\
	& = \frac{1}{N} \sum_{j \neq k} \int_{[0,T]\times \mathbb{R}^d} \bar\bfm_{t}(X^j_{[0,t]}, X^k_{[0,t]}, \clv_{[0,t]},h) \gamma(dh) dt\nonumber\\
	&\quad + \frac{1}{N} \sum_{j=1}^N \int_{[0,T]\times \mathbb{R}^d} \bar\bfm_{t}(X^j_{[0,t]}, X^j_{[0,t]}, \clv_{[0,t]},h) \gamma(dh) dt
		+   \clr^N_4, \label{eq:eqcombine}
\end{align}
where  $\clr^N_4$ 
converges to $0$ in probability as $N\to \infty$.
\end{lemma}
{\bf Proof.}  From \eqref{eq:eq628}
$$(\bfd^{N,i}_s(h) - \bfd^{i}_s(h))^2  = (\clt^{N,i}_1(s) + \clt^{N,i}_2(s) + \clt^{N,i}_3(s))^2
= \sum_{m=1}^3 (\clt^{N,i}_m(s))^2 + 2 \sum_{m < l} \clt^{N,i}_m(s) \clt^{N,i}_l(s),$$
where
\begin{align*}
\clt^{N,i}_1(s) &=  (Y^N_s-Y_s)\cdot d_2(R^i_s, h, \mu_s), \;  \clt^{N,i}_2(s) = 
 \langle  d_3(R^i_s, h, \mu_s, \cdot), (\mu^N_s-\mu_s) \rangle , \\
 \clt^{N,i}_3(s) &= \theta_d(R^i_s,  R^{N,i}_s, \mu_s, \mu^N_s, h).
\end{align*}
As for \eqref{eq:eq638} we have, as $N\to \infty$,
\begin{equation}
	\sum_{i=1}^N \int_{\mathbb{X}_T} 1_{[0,\bfd^{i}_{s-}(h)]}(u) \frac{1}{(\bfd^{i}_{s-}(h))^2}
	(\clt^{N,i}_3(s-))^2 d \bfn^i \to 0 \mbox{ in probability. }
\label{eq:eq753}	
\end{equation}
Similar estimates show that for $m = 1,2$, as $N\to \infty$,
\begin{equation}
	\sum_{i=1}^N \int_{\mathbb{X}_T} 1_{[0,\bfd^{i}_{s-}(h)]}(u) \frac{1}{(\bfd^{i}_{s-}(h))^2}
	|\clt^{N,i}_m(s-)| |\clt^{N,i}_3(s-)|  d \bfn^i \to 0 \mbox{ in probability. }
\label{eq:eq754}	
\end{equation}
Next
$$
(\clt^{N,i}_1(s))^2 = (\clt^{N,i}_4(s) + \clt^{N,i}_5(s))^2 = (\clt^{N,i}_4(s))^2 + (\clt^{N,i}_5(s))^2 + 2\clt^{N,i}_4(s)\clt^{N,i}_5(s),$$ 
where
$$\clt^{N,i}_4(s) = \frac{1}{N} \sum_{j=1}^N \bfh^j_s \cdot d_2(R^i_s, h, \mu_s), \;
 \clt^{N,i}_5(s) =
\clt^N_s \cdot d_2(R^i_s, h, \mu_s).$$
As for \eqref{eq:eq704}, we see that, as $N\to \infty$,
\begin{equation}
	\label{eq:eq804}
		\sum_{i=1}^N \int_{\mathbb{X}_T} 1_{[0,\bfd^{i}_{s-}(h)]}(u) \frac{1}{(\bfd^{i}_{s-}(h))^2}
	(	\clt^{N,i}_5(s-))^2 d \bfn^i \to 0 \mbox{ in probability. }
\end{equation}
Next, using the observation that
$\E_{P^N}\bfs_{0,t}(  X^j_{[0,t]},  v_{[0,t]})=0$, for all $j \in \N$ and $v_{[0,t]}$
in $ \DD_{\RR^{2m} \times \clp(\RR^d)}[0,t]$; and making use of Lemma \ref{lem:etbds} once more, we see that, as
$N\to \infty$,
\begin{equation}
	\label{eq:eq804b}
		\int_{\mathbb{X}_T} 1_{[0,\bfd^{i}_{s-}(h)]}(u) \frac{1}{(\bfd^{i}_{s-}(h))^2}
		|\clt^{N,i}_4(s-)\clt^{N,i}_5(s-)| d \bfn^i \to 0 \mbox{ in probability. }
		\end{equation}
Also, 
\begin{align*}
	&\sum_{i=1}^N \int_{\mathbb{X}_T} 1_{[0,\bfd^{i}_{s-}(h)]}(u) \frac{1}{(\bfd^{i}_{s-}(h))^2}
(	\clt^{N,i}_4(s-))^2 d \bfn^i	\\
 &= 	\frac{1}{N^2}\sum_{i,j,k} \int_{\mathbb{X}_T} 1_{[0,\bfd^{i}_{s-}(h)]}(u) \frac{1}{(\bfd^{i}_{s-}(h))^2}
(\bfh^j_s \cdot d_2(R^i_{s-}, h, \mu_{s-}))(\bfh^k_s \cdot d_2(R^i_{s-}, h, \mu_{s-})) d \bfn^i \\
\end{align*}
The above can be rewritten as
\begin{align}
& \frac{1}{N^2}\sum_{i,j,k} \Big(\int_{\mathbb{X}_T} 1_{[0,\bfd^{i}_{s-}(h)]}(u) \frac{1}{(\bfd^{i}_{s-}(h))^2}
(\bfh^j_s \cdot d_2(R^i_{s-}, h, \mu_{s-}))(\bfh^k_s \cdot d_2(R^i_{s-}, h, \mu_{s-})) d \bfn^i \nonumber \\
&\quad \quad \quad - \int_{[0,T]\times \mathbb{R}^d} \bar\bfm_{1,t}(X^{j}_{[0,t]}, X^{k}_{[0,t]}, \clv_{[0,t]},h) \gamma(dh) dt  \Big)\nonumber\\
&\quad + \frac{1}{N}\sum_{j\neq k} \int_{[0,T]\times \mathbb{R}^d} \bar\bfm_{1,t}(X^{j}_{[0,t]}, X^{k}_{[0,t]}, \clv_{[0,t]},h) \gamma(dh) dt 
\nonumber \\
&\quad 
+ \frac{1}{N}\sum_{j} \int_{[0,T]\times \mathbb{R}^d} \bar\bfm_{1,t}(X^{j}_{[0,t]}, X^{j}_{[0,t]}, \clv_{[0,t]},h) \gamma(dh) dt.
\label{eq:eq1100}
\end{align}
A similar argument as below \eqref{eq:eq1153} shows that the first term in the above display converges 
to $0$ in probability, as $N\to \infty$.

Combining  \eqref{eq:eq804}, \eqref{eq:eq804b} and \eqref{eq:eq1100} we have that
\begin{align}
	&\sum_{i=1}^N \int_{\mathbb{X}_T} 1_{[0,\bfd^{i}_{s-}(h)]}(u) \frac{1}{(\bfd^{i}_{s-}(h))^2}
(	\clt^{N,i}_1(s-))^2 d \bfn^i \nonumber\\
&= 	\frac{1}{N}\sum_{j\neq k} \int_{[0,T]\times \mathbb{R}^d} \bar\bfm_{1,t}(X^{j}_{[0,t]}, X^{k}_{[0,t]}, \clv_{[0,t]},h) \gamma(dh) dt \nonumber \\
	&\quad + \frac{1}{N}\sum_{j} \int_{[0,T]\times \mathbb{R}^d} \bar\bfm_{1,t}(X^{j}_{[0,t]}, X^{j}_{[0,t]}, \clv_{[0,t]},h) \gamma(dh) dt
	+ \tilde R^N_1,
	\label{eq:eq1137}
\end{align}
where $\tilde R^N_1$ 
converges to $0$ in probability as $N\to \infty$.

We now consider the term $\clt^{N,i}_2(s)$.
Writing
$$(\clt^{N,i}_2(s))^2 = \frac{1}{N^2} \sum_{j,k} d_3^c(R^i_s, h, \mu_s, X^j_s) d_3^c(R^i_s, h, \mu_s, X^k_s)$$
we see
\begin{align*}
	&\sum_{i=1}^N \int_{\mathbb{X}_T} 1_{[0,\bfd^{i}_{s-}(h)]}(u) \frac{1}{(\bfd^{i}_{s-}(h))^2}(\clt^{N,i}_2(s-))^2 d \bfn^i	\\
	&= \frac{1}{N^2}\sum_{i,j,k} \int_{\mathbb{X}_T} 1_{[0,\bfd^{i}_{s-}(h)]}(u) \frac{1}{(\bfd^{i}_{s-}(h))^2}
	d_3^c(R^i_{s-}, h, \mu_{s-}, X^j_{s-}) d_3^c(R^i_{s-}, h, \mu_{s-}, X^k_{s-}) d \bfn^i \\\
\end{align*}
The above can be rewritten as
\begin{align}
	&\frac{1}{N^2}\sum_{i,j,k} \Big(\int_{\mathbb{X}_T} 1_{[0,\bfd^{i}_{s-}(h)]}(u) \frac{1}{(\bfd^{i}_{s-}(h))^2}
	d_3^c(R^i_{s-}, h, \mu_{s-}, X^j_{s-}) d_3^c(R^i_{s-}, h, \mu_{s-}, X^k_{s-}) d \bfn^i \nonumber \\
	&\quad \quad \quad - \int_{[0,T]\times \mathbb{R}^d} \bar\bfm_{3}(X^{j}_{t}, X^{k}_{t}, \clv_{t},h) \gamma(dh) dt  \Big)\nonumber\\
	&\quad + \frac{1}{N}\sum_{j\neq k} \int_{[0,T]\times \mathbb{R}^d} \bar\bfm_{3}(X^{j}_{t}, X^{k}_{t}, \clv_{t},h) \gamma(dh) dt 
	+ \frac{1}{N}\sum_{j} \int_{[0,T]\times \mathbb{R}^d} \bar\bfm_{3}(X^{j}_{t}, X^{j}_{t}, \clv_{t},h) \gamma(dh) dt.	
	\label{eq:eq1241}
\end{align}
As before, the first term above 
converges to $0$ in probability, as $N\to \infty$.
 Thus
\begin{align}
	&\sum_{i=1}^N \int_{\mathbb{X}_T} 1_{[0,\bfd^{i}_{s-}(h)]}(u) \frac{1}{(\bfd^{i}_{s-}(h))^2}(\clt^{N,i}_2(s-))^2 d \bfn^i	\nonumber\\
	&\quad = \frac{1}{N}\sum_{j\neq k} \int_{[0,T]\times \mathbb{R}^d} \bar\bfm_{3}(X^{j}_{t}, X^{k}_{t}, Y_{t}, \mu_{t},h) \gamma(dh) dt \nonumber \\
	&\quad + \frac{1}{N}\sum_{j} \int_{[0,T]\times \mathbb{R}^d} \bar\bfm_{3}(X^{j}_{t}, X^{j}_{t}, Y_{t}, \mu_{t},h) \gamma(dh) dt + \tilde R^N_2,
\label{eq:eq1137b}	
\end{align}
where $\tilde R^N_2$ 
converges to $0$ in probability as $N\to \infty$.

We now consider the term $2\clt^{N,i}_1(s)\clt^{N,i}_2(s)$.
\begin{align}
2\clt^{N,i}_1(s)\clt^{N,i}_2(s) &= 	2(Y^N_s-Y_s)\cdot d_2(R^i_s, h, \mu_s)
	 \langle  d_3(R^i_s, h, \mu_s, \cdot), (\mu^N_s-\mu_s) \rangle	\nonumber\\
	&= \frac{2}{N^2}\sum_{j,k} \bfh^j_s \cdot d_2(R^i_s, h, \mu_s)d_3^c(R^i_s, h, \mu_s, X^k_s)\nonumber\\
	&\quad + \frac{2}{N}\sum_{k}\clt^N_s\cdot d_2(R^i_s, h, \mu_s) 
	d_3^c(R^i_s, h, \mu_s, X^k_s)\nonumber\\
	&\quad \equiv \clt^{N,i}_6(s) + \clt^{N,i}_7(s).
	\label{eq:eq1255}
\end{align}
For the  term $\clt^{N,i}_6(s)$ note that,
\begin{align}
	&\sum_{i=1}^N \int_{\mathbb{X}_T} 1_{[0,\bfd^{i}_{s-}(h)]}(u) \frac{1}{(\bfd^{i}_{s-}(h))^2}\clt^{N,i}_6(s-) d \bfn^i	\nonumber\\
	&\quad = \frac{2}{N^2}\sum_{i,j,k} \int_{\mathbb{X}_T} 1_{[0,\bfd^{i}_{s-}(h)]}(u) \frac{1}{(\bfd^{i}_{s-}(h))^2}
	\bfh^j_s \cdot d_2(R^i_{s-}, h, \mu_{s-})d_3^c(R^i_{s-}, h, \mu_{s-}, X^k_{s-}) d \bfn^i
\end{align}
As in \eqref{eq:eq1241} and \eqref{eq:eq1137b}, we can now write the above as
\begin{align}
	&\sum_{i=1}^N \int_{\mathbb{X}_T} 1_{[0,\bfd^{i}_{s-}(h)]}(u) \frac{1}{(\bfd^{i}_{s-}(h))^2}\clt^{N,i}_6(s-) d \bfn^i	\nonumber\\
	&\quad = \frac{1}{N}\sum_{j\neq k} \int_{[0,T]\times \mathbb{R}^d} \bar\bfm_{2,t}(X^{j}_{[0,t]}, X^{k}_{[0,t]}, Y_{[0,t]}, \mu_{[0,t]},h) \gamma(dh) dt \nonumber \\
	&\quad + \frac{1}{N}\sum_{j} \int_{[0,T]\times \mathbb{R}^d} \bar\bfm_{2,t}(X^{j}_{[0,t]}, X^{j}_{[0,t]}, Y_{[0,t]}, \mu_{[0,t]},h) \gamma(dh) dt + \tilde R^N_3,
\label{eq:eq1137c}	
\end{align}
where  $\tilde R^N_3$ 
converges to $0$ in probability as $N\to \infty$.
Also, as for \eqref{eq:eq804b}, as $N\to \infty$,
$$\sum_{i=1}^N \int_{\mathbb{X}_T} 1_{[0,\bfd^{i}_{s-}(h)]}(u) \frac{1}{(\bfd^{i}_{s-}(h))^2}|\clt^{N,i}_7(s-)| d \bfn^i	 \to 0 \mbox{ in probability. }$$
The result now follows on combining the above display with \eqref{eq:eq753}, \eqref{eq:eq754}, \eqref{eq:eq1137}, \eqref{eq:eq1137b}, \eqref{eq:eq1255}, 
\eqref{eq:eq1137c}.	
\qed

Recall the function $\vartheta^{N,i}_{s}$ introduced at the beginning of the subsection. Using very similar estimates as in the proof of Lemma \ref{lem:lem910}, one can establish the following result.  We omit the proof.
\begin{lemma}\label{lem:lem5120}
As $N\to \infty$,
	\begin{align*}
	\sum_{i=1}^N \int_{\mathbb{X}_T} 1_{[0,\bfd^{i}_{s-}(h)]}(u) \vartheta^{N,i}_{s-}(h)\left(\frac{\bfd^{N,i}_{s-}(h)}{\bfd^{i}_{s-}(h)}-1\right)^3 d \bfn^i
\end{align*}
converges to $0$ in probability.
\end{lemma}
\subsection{Comment on Smoothness Conditions.}
\label{sec:remoncond}
Conditions \ref{cond:clt2}, \ref{cond:clt3} and \ref{cond:clt4} on $b_0$, $b$ and $d$ can be regarded as smoothness  conditions.
These conditions  are satisfied quite generally. We  give two examples to illustrate this.
\begin{example}
	\label{ex:ex747}
	Let $d= m =1$. Let $\bar{b}:\mathbb{R}^{k+2}\rightarrow \mathbb{R}$ be bounded Lipschitz and twice continuously differentiable, with bounded derivatives, in the last $k+1$ variables. 
	 Let $\bar{b}_0:\mathbb{R}^{k+1}\rightarrow \mathbb{R}$ be bounded Lipschitz and twice continuously differentiable with bounded  derivatives. Similar assumptions on $\bar{\sigma}_0^i$ for $i=1,\cdots,m$. Let $\bar{d}:\mathbb{R}^{k+3} \to (\epsilon, \infty)$ be bounded and Lipschitz in the first $k+2$ variables, uniformly in the last variable, where $\epsilon \in (0, \infty)$.   Also suppose that
	$\bar d$ is  twice continuously differentiable, with bounded derivatives, with respect to the middle $k+1$ variables.  Now let $b$, $b_0$ and $d$ be of the form: For $(x,y,\nu,h) \in \RR^{d+m}\times \clp(\RR^d)\times \RR^d$
\begin{itemize}
\item $b(x,y,\nu)=\bar{b}(x,y,\langle f_1,\nu\rangle,\cdots,\langle f_k,\nu\rangle)$,
\item $b_0(y,\nu)=\bar{b}_0(y,\langle f_1,\nu\rangle,\cdots,\langle f_k,\nu\rangle)$,
\item $\sigma_0^i(y,\nu)=\bar{\sigma}^i_0(y,\langle f_1,\nu\rangle,\cdots,\langle f_k,\nu\rangle)$,
\item $d(x,y,\nu,h)=\bar{d}(x,y,\langle f_1,\nu\rangle,\cdots,\langle f_k,\nu\rangle,h)$,
\end{itemize}where $f_i$ are bounded Lipschitz functions. Finally let 
$d_0: \RR \to (\epsilon, \infty)$ be a bounded function and let
$\gamma$, $\gamma_0$ be  probability measures on $\mathbb{R}$ with finite second moment. Then it is easy to check that Conditions \ref{cond:main} and \ref{cond:clt1}  is satisfied. For Condition \ref{cond:clt2} observe that by Taylor's expansion,
\begin{align*}
	b(z',\nu')-b(z,\nu) &=(y'-y)\bar{b}_y(z,\langle f_1,\nu\rangle,\cdots,\langle f_k,\nu\rangle)\\
	&\quad + \sum_{i=1}^k\bar{b}_{u_i}(z,\langle f_1,\nu\rangle,\cdots,\langle f_k,\nu\rangle) \langle f_i,(\nu'-\nu)\rangle + \theta_b(z,z',\nu,\nu')\,,
\end{align*}
where for some constant $K_1$, 
$|\theta_b(z,z',\nu,\nu')|\leq K_1 (|y'-y|^2 + \displaystyle\max_{i}| \langle f_i,(\nu'-\nu)\rangle|^2)$. 
This verifies  Condition \ref{cond:clt3}.  Conditions \ref{cond:clt2}, \ref{cond:clt4} can  be verified similarly.
\end{example}
\begin{example}Let $d= m =1$. Let $\tilde{b}:\mathbb{R}^3\rightarrow \mathbb{R}$, $\tilde{b}_0:\mathbb{R}^2\rightarrow \mathbb{R}$, $\tilde{\sigma}^i_0:\mathbb{R}^2\rightarrow \mathbb{R}, i=1,\cdots,m$ be bounded Lipschitz functions. Further suppose that $\tilde{b}$ is twice continuously differentiable with respect to the second variable with bounded  derivatives and $\tilde{b}_0$ is also twice continuously differentiable with respect to the first variable,  with bounded  derivatives. Similar assumptions on $\tilde{\sigma}^i_0$. Let $\tilde{d}:\mathbb{R}^4\rightarrow \mathbb{R}_+$ be bounded and Lipschitz in the first three variables, uniformly in the last variable. 
	Also suppose that $\tilde d$ is twice continuously differentiable, with bounded  derivatives,  in the second variable.
	Let $d_0, \gamma, \gamma_0$ be as in Example \ref{ex:ex747}.
Now let $b$, $b_0$ and $d$ be of the form:
\begin{itemize}
\item $b(x,y,\nu)=\int\tilde{b}(x,y,x')\nu(dx')$,
\item $b_0(y,\nu)=\int\tilde{b}_0(y,x')\nu(dx')$,
\item $\sigma^i_0(y,\nu)=\int\tilde{\sigma}^i_0(y,x')\nu(dx')$,
\item $d(x,y,\nu,h)=\int\tilde{d}(x,y,x',h)\nu(dx')$.
\end{itemize}
 Then it is easy to check that for this example Condition \ref{cond:main} is satisfied. One can also check that Conditions \ref{cond:clt2}, \ref{cond:clt3} and \ref{cond:clt4}  are satisfied as well.  In particular, note that for $x\in \RR^d$, $y, y' \in \RR^m$, $\nu, \nu' \in \clp(\RR^m)$,
\begin{align*}
b(x,y',\nu')-b(x,y,\nu)&=\int(\tilde{b}(x,y,x')(\nu'-\nu)(dx')
+\int(\tilde{b}(x,y',x')-\tilde{b}(x,y,x'))\nu(dx')
\\&+\int(\tilde{b}(x,y',x')-\tilde{b}(x,y,x'))(\nu'-\nu)(dx')\,.
\end{align*}
Using Taylor's expansion to the second term we get,
\begin{align*}
\int(\tilde{b}(x,y',x')-\tilde{b}(x,y,x'))\nu(dx')= (y'-y) \int\tilde{b}_y(x,y,x')\nu(dx')+\frac{1}{2}(y'-y)^2
r_1(x,y,y'),
\end{align*}
where $r_1$ is a bounded function.
Using Taylor's expansion to the third term we get
\begin{align*}
\int(\tilde{b}(x,y',x')-\tilde{b}(x,y,x'))(\nu'-\nu)(dx')&= (y'-y)\int\tilde{b}_y(x,y,x')(\nu'-\nu)(dx')\\
&+ \frac{1}{2}(y'-y)^2 r_2(x,y,y'),
\end{align*}
where $r_2$ is a bounded function.
Finally using the boundedness and continuity of $\tilde b$, $\tilde b_y$ and the  inequality
$$|(y'-y)\int\tilde{b}_y(x,y,x')(\nu'-\nu)(dx')|\leq |y'-y|^2+|\int\tilde{b}_y(x,y,x')(\nu'-\nu)(dx')|^2
$$ we see that Condition \ref{cond:clt3} is satisfied.  Conditions \ref{cond:clt2},\ref{cond:clt4}  can be verified similarly.
\end{example}
\subsection{Some Integral Operators.}
\label{sec:someintegop}
Define for $t \in [0,T]$, the function $\bff_t$ from $\DD_{\RR^{2d+2m} \times \clp(\RR^d)}[0,t]$ to $\RR^d$ as follows:\\  
For $(x^{(1)}_{[0,t]}, x^{(2)}_{[0,t]}, y_{[0,t]}, w_{[0,t]}, \nu_{[0,t]})= \zeta_{[0,t]} \in  \DD_{\RR^{2d+m} \times \clp(\RR^d)}[0,t]$ 
\begin{equation}\bff_{t}(\zeta_{[0,t]}) =
b_{3}^c(x^{(1)}_{t},y_t, \nu_t, x^{(2)}_{t}) + 
b_2(x^{(1)}_{t}, y_t, \nu_t)\bfs_{0,t}(\zeta^{(2)}_{[0,t]}),\label{eq:bfft}\end{equation}
where as before $\zeta^{(2)} = (x^{(2)}, y, w, \nu)$.
We note that
\begin{align}\|\bff_{t}(\zeta_{[0,t]})\|^2 &=
	\bfs_{1,t}(x^{(1)}_{t}, x^{(2)}_{[0,t]},  \zeta^{(2)}_{[0,t]})
	+ \bfs_{2,t}(x^{(1)}_{t}, x^{(2)}_{[0,t]}, \zeta^{(2)}_{[0,t]})
	+ \bfs_{3}(x^{(1)}_{t}, x^{(2)}_{t}, x^{(2)}_{t}, y_{[0,t]}, \nu_{[0,t]}).
	\label{eq:eq219new}
\end{align}
  Also define for $t \in [0,T]$, the function $\bar\bff_t$ from $\DD_{\RR^{2d+2m} \times \clp(\RR^d)}[0,t]\times \RR_+ \times \RR^d$ to $\RR$ as follows:
For $(x^{(1)}_{[0,t]}, x^{(2)}_{[0,t]},y_{[0,t]},w_{[0,t]}, \nu_{[0,t]},u,h)= (\zeta_{[0,t]},u,h) \in  \DD_{\RR^{2d+m} \times \clp(\RR^d)}[0,t]\times \RR_+ \times \RR^d$
\begin{align*}
	\bar\bff_{t}(\zeta_{[0,t]},u,h) &=
	1_{[0,d(x^{(1)}_{t-},y_{t-},\nu_{t-},h)]}(u)\frac{1}{d(x^{(1)}_{t-},y_{t-},\nu_{t-},h)}\Big (d_3^c(x^{(1)}_{t-},y_{t-},h,\nu_{t-},x^{(2)}_{t-})\\
	&\quad +  \bfs_{0,t}(\zeta^{(2)}_{[0,t]}) \cdot d_2(x^{(1)}_{t-},y_{t-},h,\nu_{t-})\Big ).
\end{align*}
The functions $\bff_t, \bar \bff_t$ will play the role of kernels for certain integral operators on $L^2$ spaces.
To describe these operators,
in addition to the canonical  spaces and processes introduced in Section \ref{sec:cltassu} (see \eqref{eq:canon1}, \eqref{eq:canon2}), we define the canonical processes
 $V_*^0 = (B_*^0, \bfn_*^0, Y_*)$ on  $\Om_m$  as
\begin{align*}
V_*^0(\om_0) &= (B_*^0(\om_0), \bfn_*^0(\om_0), Y_*(\om_0)) = (\om_{0,1},\om_{0,2},\om_{0,3});\; \om_0 = (\om_{0,1},\om_{0,2},\om_{0,3}) \in \Om_m.
\end{align*}
Also, with $\Pi$ as introduced in Remark \ref{rem:rem415}, let $\mu_*: \Om_m \to \mathbb{D}_{\clp(\RR^d)}[0,T]$ be defined as
$\mu_*(\om_0) = \Pi(Y_{*0}(\om_0),B_*^0(\om_0), \bfn_*^0(\om_0))$.
Write
\begin{equation}
	\clv_* = (Y_*, B_*^0, \mu_*). \label{eq:clvstar}
\end{equation}
We can now define the integral operators related to $\bff_t$ and  $\bar \bff_t$. 
Recall the transition probability kernel $\alpha$ introduced in \eqref{eq:transkern}.
Fix $\om_0 \in \Om_m$ and consider the Hilbert space $\clh_{\om_0} = L^2(\Om_d, \alpha(\om_0, \cdot))$.
We denote the norm and inner product in $\clh_{\om_0}$ as $\|\cdot\|_{\om_0}$ and $\langle \cdot, \cdot\rangle_{\om_0}$ respectively.
Define the integral operator $A^1_{\om_0}$ on $\clh_{\om_0}$ as follows.  For $g^{(1)} \in \clh_{\om_0}$, $(A^1_{\om_0} g^{(1)}) = \hat g^{(1)}_{\om_0}$, where for $\om_2 \in \Om_d$,
$$
\hat g^{(1)}_{\om_0}(\om_2) = \int_{\Om_d} g^{(1)}(\om_1) \left( \int_0^T \bff_t(X_{*,[0,t]}(\om_1), X_{*,[0,t]}(\om_2), \clv_{*,[0,t]}(\om_0)) dB_{*,t}(\om_1) \right) \alpha(\om_0, d\om_1).$$
Also define the integral operator $A^2_{\om_0}$ on $\clh_{\om_0}$ as follows.  For $g^{(2)} \in \clh_{\om_0}$, $(A^2_{\om_0} g^{(2)}) = \hat g^{(2)}_{\om_0}$, where for $\om_2 \in \Om_d$,
$$
\hat g^{(2)}_{\om_0}(\om_2) = \int_{\Om_d} g^{(2)}(\om_1)\left(\int_{\XX_T} \bar \bff_t(X_{*,[0,t]}(\om_1), X_{*,[0,t]}(\om_2), \clv_{*,[0,t]}(\om_0),u,h) d\tilde \bfn_*(\om_1) \right) \alpha(\om_0, d\om_1).$$
Let $A_{\om_0} = A^1_{\om_0} + A^2_{\om_0}$.
Denote by $I$ the identity operator on $\clh_{\om_0}$.
\begin{lemma}
	\label{lem:lem916}
	For $P_0$ a.e. $\om_0$, (i) $\mbox{Trace}(A^1_{\om_0}(A^2_{\om_0})^*) = 0$; (ii) $\mbox{Trace}(A^n_{\om_0}) = 0$ for all $n \ge 2$; and
(iii) $I- A_{\om_0}$ is invertible.
\end{lemma}
{\bf Proof.} 
Parts (i) and (ii) are consequences of independence between  $B_*$ and $\bfn_*$ under $\alpha(\om_0,\cdot)$.
For example for (i), from the definitions of $A^i_{\om_0}$, it follows that
\begin{align*}\mbox{Trace}(A^1_{\om_0}(A^2_{\om_0})^*)=\int_{\Omega_d^2} \left( \int_0^T \bff_t(X_{*,[0,t]}(\om_1), X_{*,[0,t]}(\om_2), \clv_{*,[0,t]}(\om_0)) dB_{*,s}(\om_1) \right)\\ \left(\int_{\XX_T} \bar \bff_t(X_{*,[0,t]}(\om_1), X_{*,[0,t]}(\om_2), \clv_{*,[0,t]}(\om_0),u,h) d\tilde \bfn_*(\om_1) \right)\alpha(\om_0,d\om_1)\alpha(\om_0,d\om_2)\,.
\end{align*}
The above expression is $0$ due to the independence between  $B_*$ and $\bfn_*$ under $\alpha(\om_0,\cdot)$.
Part (ii) is proved similarly (see e.g. Lemma 2.7 of \cite{ShiTan} ).
Part  (iii)  is now immediate from Lemma 1.3 of \cite{ShiTan}.
\qed

\subsection{Combining Contributions from $J^{N,1}$ and $J^{N,2}$.}
Recall the integral operators $A^i_{\om_0}$, $i=1,2$, introduced in Section \ref{sec:someintegop}.
Define $\tau^{(i)}:  \Om_m \to \RR$ as $\tau^{(i)}(\om_0) =  \mbox{Trace}(A^i_{\om_0}(A^i_{\om_0})^*)$, $i=1,2$.
From Lemma \ref{lem:lem916}  we have that, for $P_0$ a.e. $\om_0$,
\begin{equation}
	\mbox{Trace}(A_{\om_0}(A_{\om_0})^*) = \tau^{(1)}(\om_0) + \tau^{(2)}(\om_0).
\end{equation}
The following lemma gives the asymptotics for the second terms on the right sides of \eqref{eq:eqquadform}
and \eqref{eq:eqcombine}.
\begin{lemma}
	\label{lem:lemtrace}
	As $N \to \infty$, 
	$$ \frac{1}{N} \sum_{j=1}^N \int_0^T \bfm_{t}(X^j_{[0,t]}, X^j_{[0,t]}, \clv_{[0,t]}) dt -\tau^{(1)}(V^0)
	$$
	and
	$$ \frac{1}{N} \sum_{j=1}^N \int_{[0,T]\times \mathbb{R}^d} \bar\bfm_{t}(X^j_{[0,t]}, X^j_{[0,t]}, \clv_{[0,t]},h)\gamma(dh) dt - \tau^{(2)} (V^0)
	$$
converge to $0$ in probability.
\end{lemma}
{\bf Proof.} Note that if $A$ is an integral operator on $L^2(\nu)$ with associated kernel $a(x,y)$, then $\mbox{Trace}(AA^*)=||a||^2_{L^2(\nu\otimes \nu)}$. Thus from the definition of the operator $A^1_{\om_0}$, 
\begin{align*}
	&\mbox{Trace}(A^1_{\om_0}(A^1_{\om_0})^*)\\
	&\quad =\int_{\Omega_d^2}\bigg|\int_0^T \bff_t(X_{*,[0,t]}(\om_1), X_{*,[0,t]}(\om_2), \clv_{*,[0,t]}(\om_0)) dB_{*,t}(\om_1)\bigg|^2\alpha(\om_0,d \om_1)\alpha(\om_0,d\om_2)\\
	&\quad =\int_{\Omega_d}\int_0^T\int_{\Omega_d}\|\bff_t(X_{*,[0,t]}(\om_1), X_{*,[0,t]}(\om_2), \clv_{*,[0,t]}(\om_0))\|^2\alpha(\om_0,d\om_1)\, dt\, \alpha(\om_0,d\om_2)\\
\end{align*}
Using the relation \eqref{eq:eq219new} we have,
\begin{align*}
	\mbox{Trace}(A^1_{\om_0}(A^1_{\om_0})^*)
	& =\int_{\Omega_d}\int_0^T\int_{\Omega_d}\Big(\bfs_{1,t}(X_{*,t}(\om_1),  X_{*,[0,t]}(\om_2), X_{*,[0,t]}(\om_2), \clv_{*,[0,t]}(\om_0))\\
	&\quad +\bfs_{2,t}(X_{*,t}(\om_1),  X_{*,[0,t]}(\om_2), X_{*,[0,t]}(\om_2), \clv_{*,[0,t]}(\om_0))\\
	&\quad +\bfs_{3}(X_{*,t}(\om_1),  X_{*,t}(\om_2), X_{*,t}(\om_2), \clv_{*,t}(\om_0))\Big)\alpha(\om_0,d\om_1)\, dt\, \alpha(\om_0,d\om_2)\\
	&  =
\int_{\Omega_d} \int_0^T \bfm_{t}(X_{*,[0,t]}(\om_2),X_{*,[0,t]}(\om_2), \clv_{*,[0,t]}(\om_0))\, dt\, \alpha(\om_0,d \om_2).
\end{align*}

Since conditional on $\clg^0$, $\{X^j\}$ are i.i.d. with common distribution
$\alpha(V^0, \cdot)\circ X_*^{-1}$, the first convergence in the lemma now follows from the weak law of large numbers.  The second convergence statement
is proved similarly.
\qed


We will now use the results from Section \ref{sec:symmstat} with $\clx = \Om_d$ and $\nu = \alpha(\om_0, \cdot)$, $\om_0 \in \Om_m$.
For each $\om_0 \in \Om_m$, $k \ge 1$ and $f \in L^2_{sym}(\alpha(\om_0, \cdot)^{\otimes k})$ the multiple stochastic integral $I_k^{\om_0}(f)$
is defined as in Section \ref{sec:symmstat}.  More precisely, let
$\cla^p$ be the collection of all measurable $f: \Om_{m}\times \Om_d^p \to \RR$ such that
$$\int_{\Om^p}|f(\om_0, \om_1, \ldots , \om_p)|^2 \alpha(\om_0, d\om_1)\cdots \alpha(\om_0, d\om_p) < \infty, \; P_0 \mbox{ a.e. } \om_0$$
and $f(\om_0, \cdot)$ is symmetric for $P_0$  a.e. $\om_0$.
Then there is a measurable space $(\Om^*, \clf^*)$ and a regular conditional probability distribution $\alpha^*: \Om_0 \times \clf^* \to [0,1]$
such that on the probability space 
$
(\Om_m \times \Om^*, \clb(\Om_m)\otimes \clf^*, P_0 \otimes \alpha^*)$,
where
$$P_0 \otimes \alpha^*(A \times B) = \int_A \alpha^*(\om_0, B) P_0(d\om_0), \; A\times B \in \clb(\Om_m)\otimes \clf^*,
$$
there is a collection or real valued random variables $\{I_p(f): f \in \cla^p, p \ge 1\}$ with the properties that
\begin{enumerate}[(a)]
	\item For all $f \in \cla^1$ the conditional distribution of $I_1(f)$ given $\clg^0_* = \clb(\Om_m)\otimes \{\emptyset,\Om^*\}$
	is Normal with mean $0$ and variance $\int_{\Om_d} f^2(\om_0, \om_1) \alpha(\om_0, d\om_1)$.
	\item $I_p$ is (a.s.) linear map on $\cla^p$.
	\item For $f \in \cla^p$ of the form
	$$f(\om_0, \om_1, \ldots , \om_p) = \prod_{i=1}^p h(\om_0, \om_i), \mbox{ s.t. }  \int_{\Om_d} h^2(\om_0, \om_1) \alpha(\om_0, d\om_1) < \infty, \; P_0 \mbox{ a.e. } \om_0,$$
	$$I_p(f) = \sum_{j=0}^{\lfloor p/2 \rfloor} (-1)^j C_{p,j} \left( \int_{\Om_d} h^2(\om_0, \om_1) \alpha(\om_0, d\om_1)\right)^j (I_1(h))^{p-2j}$$
 and
	$$\int_{\Om^*} \left(I_p(f)(\om_0, \om^*)\right)^2 \alpha^*(\om_0, d\om^*)
	= p! \left( \int_{\Om_d} h^2(\om_0, \om_1) \alpha(\om_0, d\om_1)\right)^p$$
\end{enumerate}
$P_0$ a.e.  $\om_0$.
We write $I_p(f)(\om_0, \cdot)$ as $I_p^{\om_0}(f)$.
With an abuse of notation, we will denote once more by $V^0_*$ the canonical process on $\Om_m \times \Om^*$, i.e. $V^0_*(\om_0, \om^*) = \om_0$, for
$(\om_0, \om^*) \in \Om_m \times \Om^*$.

  Recall the class $\cla$ introduced in Section \ref{sec:cltassu}.  Let for $\varphi \in \cla$
$$
\clv_N^{\varphi} = \sqrt{N} \left (\frac{1}{N} \sum_{j=1}^N \varphi(X^{j}) - m_{\varphi}(V^0)\right).$$
Define $\bar \tau: \Om_m \times \Om^* \to \RR$ as
$\bar \tau(\om_0, \om^*) =  \mbox{Trace}(A_{\om_0}(A_{\om_0})^*)$.

Given $\om_0 \in \Om_m$, define
$F^1_{\om_0}: \Om_d \times \Om_d \to \RR$ as follows: For $(\om_1, \om_2) \in \Om_d \times \Om_d$
\begin{align*}
	F^1_{\om_0}(\om_1, \om_2) &=
	\int_0^T \bff_t(X_{*,[0,t]}(\om_1), X_{*,[0,t]}(\om_2), \clv_{*,[0,t]}(\om_0)) dB_{*,t}(\om_1)\\
	&\quad + \int_0^T \bff_t(X_{*,[0,t]}(\om_2), X_{*,[0,t]}(\om_1), \clv_{*,[0,t]}(\om_0)) dB_{*,t}(\om_2)\\
	&\quad -  \int_0^T \bfm_t(X_{*,[0,t]}(\om_1), X_{*,[0,t]}(\om_2), \clv_{*,[0,t]}(\om_0)) dt.
\end{align*}

Also, given $\om_0 \in \Om_m$, define
$F^2_{\om_0}: \Om_d \times \Om_d \to \RR$ as follows: For $(\om_1, \om_2) \in \Om_d \times \Om_d$
\begin{align*}
	F^2_{\om_0}(\om_1, \om_2) &=
	\int_{\XX_T} \bar \bff_t(X_{*,[0,t]}(\om_1), X_{*,[0,t]}(\om_2), \clv_{*,[0,t]}(\om_0),u,h) d\tilde \bfn_*(\om_1)\\
	&\quad + \int_{\XX_T} \bar\bff_t(X_{*,[0,t]}(\om_2), X_{*,[0,t]}(\om_1), \clv_{*,[0,t]}(\om_0)) d\tilde \bfn_*(\om_2)\\
	&\quad - \int_{[0,T]\times \RR^d} \bar\bfm_t(X_{*,[0,t]}(\om_1), X_{*,[0,t]}(\om_2), \clv_{*,[0,t]}(\om_0),h) \gamma(dh) dt,
\end{align*}
where  $\tilde\bfn_*$ is the compensated PRM: $\tilde \bfn_* = \bfn_* - \bfnu$.

Also let $F: \Om_m \times \Om_d \times \Om_d \to \RR$ be defined
as
$$F(\om_0, \om_1, \om_2) = \frac{1}{2} \left(F^1_{\om_0}(\om_1, \om_2) + F^2_{\om_0}(\om_1, \om_2)\right), \; (\om_0, \om_1, \om_2) \in \Om_m \times \Om_d \times \Om_d .$$
Let
$$\sigma_2^N(F) = \sum_{\substack{i,j=1\\ i \neq j}}^N F(V^0, V^i, V^j).$$
From Lemmas \ref{lem:term1}, \ref{lem:term2}, \ref{lem:lem550}, \ref{lem:lem910} and \ref{lem:lem5120} it follows that
\begin{align}
	J^{N,1}(T) + J^{N,2}(T)
	&  =  N^{-1}\sigma_2^N(F)
	 - 
	\frac{1}{2N}\sum_{j=1}^N \int_0^T \bfm_{t}(X^j_{[0,t]}, X^j_{[0,t]}, \clv_{[0,t]}) dt \nonumber\\
	&\quad - \frac{1}{2N}\sum_{j=1}^N \int_{[0,T]\times \RR^d} \bar\bfm_{t}(X^j_{[0,t]}, X^j_{[0,t]}, \clv_{[0,t]},h)\gamma(dh) dt\Big)
	+ \mathcal{R}^N, \label{eq:eq922new}
\end{align}
where $\clr^N$ converges to $0$ in probability.

In order to study the asymptotics of the expression on the left side of \eqref{eq:eq209}, we need to consider the joint asymptotic behavior of $\clv_N^{\varphi}$ and $N^{-1}\sigma_2^N(F)$.
Denote by $\bfell_{\varphi}^N$ the measurable map from $\Om_m$ to $\clp(\RR^2)$ such that
$$\cll\left( (\clv_N^{\varphi}, N^{-1}\sigma_2^N(F)) \mid \clg^0\right) = \bfell_{\varphi}^N(V^0), \mbox{ a.s. }$$

Next note that $F \in \cla^2$ and so $I_2(F)$ is a well defined random variable 
on $(\Om_m \times \Om^*, \clb(\Om_m)\otimes \clf^*, P_0 \otimes \alpha^*)$.
Also define $\bar \Phi: \Om_m \times \Om_d \to \RR$ as
$$
\bar \Phi(\om_0, \om_1) = \Phi_{\om_0}(\om_1) = \varphi(X_*(\om_1)) - m_{\varphi}(\om_0).$$
 Note that $\bar \Phi \in \cla^1$ and so
$I_1(\bar \Phi)$ is well defined.
Let $\bfell_{\varphi}$ be a measurable map from $\Om_m$ to $\clp(\RR^2)$ such that
$$
\cll\left( (I_1(\bar \Phi), I_2(F)) \mid \clg^0_*\right) = \bfell_{\varphi}(V_*^0).$$
From Theorem \ref{thm:DM} it follows that
\begin{equation}
	\label{eq:eq618}
	\bfell_{\varphi}^N(\om_0) \to \bfell_{\varphi}(\om_0) \mbox{ weakly for } P_0 \mbox{ a.e. } \om_0.
\end{equation}
The following lemma is the key step.
\begin{lemma}
\label{lem:lemcgceexp}
As $N \to \infty$,
$i \clv_N^{\varphi} + J^{N,1}(T) + J^{N,2}(T)$	
converges in distribution to $i I_1(\bar \Phi) + \frac{1}{2}I_2(F) - \frac{1}{2} \bar \tau$.
\end{lemma}
{\bf Proof.}
From \eqref{eq:eq922new} and Lemma \ref{lem:lemtrace} we have that
\begin{align*}
	i \clv_N^{\varphi} + J^{N,1}(T) + J^{N,2}(T)
	 = i\clv_N^{\varphi} + N^{-1}\sigma_2^N(F)- \frac{1}{2} (\tau^{(1)}(V^0) + \tau^{(2)}(V^0))
	+ \tilde{\mathcal{R}}^N,
\end{align*}
where $\tilde{\clr}^N$ converges to $0$ in probability.
There are measurable maps $\bfzeta^N$, $\bfzeta$ from $\Om_m$ to $\clp(\mathbb{C})$, where $\mathbb{C}$ is the complex plane, such that
with
$$
\cls^N = i\clv_N^{\varphi} + N^{-1}\sigma_2^N(F)- \frac{1}{2} (\tau^{(1)}(V^0) + \tau^{(2)}(V^0))$$
and
$$
\cls = iI_1(\bar \Phi) + I_2(F) - \frac{1}{2} \bar \tau$$
$$
\cll(\cls^N \mid \clg^0) = \bfzeta^N(V^0), \; \cll(\cls \mid \clg^0_*) = \bfzeta(V_*^0).
$$
From \eqref{eq:eq618} and the definitions of $\tau^{(i)}$ and $\bar \tau$,
\begin{equation}
	\label{eq:eq651}
	\bfzeta^N(\om_0) \to \bfzeta(\om_0), \mbox{ weakly for } P_0 \mbox{ a.e. } \om_0. 
\end{equation}
Finally, denote the probability distribution of
$(V^0, \cls^N)$ on $\Om_m \times \mathbb{C}$ by $\rho^N$ and that of $(V_*^0, \cls)$ on $\Om_m \times \mathbb{C}$ by $\rho$.  Then $\rho^N$ and $\rho$
can be disintegrated as
$$
\rho^N(A \times B) = \int_A \bfzeta^N(\om_0)(B) P_0(d\om_0), \; \rho(A \times B) = \int_A \bfzeta(\om_0)(B) P_0(d\om_0),$$
for $A \in \clb(\Om_m)$, $B \in \clb(\mathbb{C})$.
From \eqref{eq:eq651} it now follows that $\rho^N \to \rho$ weakly.  The result follows.
\qed

\subsection{Completing the proof of Theorem \ref{thm:clt}.}
\label{sec:compprf}
Recall the operator $A_{\om_0}$ introduced in Section \ref{sec:someintegop} and let $\Phi_{\om_0}$ be as in \eqref{eq:7.a}.
Define for $\om_0 \in \Om_m$,
\begin{equation}\sigma_{\om_0}^{\varphi} = \|(I-A_{\om_0})^{-1}\Phi_{\om_0}\|_{L^2(\Omega_d,\alpha(\omega_0,\cdot))}.\label{eq:eq1209}\end{equation}
	
 It follows from Lemma 1.2 of \cite{ShiTan} and Lemma \ref{lem:lem916} that  $P_0$ a.s.
 $$\E[\exp(\frac{1}{2}I_2(F))\mid \clg^0_*] = \exp(\frac{1}{2}\mbox{Trace}(A_{V^0_*}(A_{V^0_*})^*))$$ 
where $\E$ is the expectation operator on $(\Om_m \times \Om^*, \clb(\Om_m)\otimes \clf^*, P_0 \otimes \alpha^*)$.
Therefore
$$
\E\exp\left(\frac{1}{2}I_2(F) - \frac{1}{2}\bar \tau \right) = 1.$$
Also, recall that
$$
\E_{P^N}\exp\left(J^{N,1}(T) + J^{N,2}(T)\right) = 1.$$
Now applying Lemma \ref{lem:lemcgceexp} with $\varphi \equiv 0$ and using the above two displays along with 
 Scheffe's theorem we have that
$\exp(J^{N,1}(T) + J^{N,2}(T))$ is uniformly integrable.  Also since $|\exp(i\clv_N^{\varphi})|=1$,
$$\exp(i\clv_N^{\varphi} + J^{N,1}(T) + J^{N,2}(T))$$
is uniformly integrable as well.  Using Lemma \ref{lem:lemcgceexp} again we have that
 \begin{align*}
 &\lim_{N \rightarrow \infty}\E_{P^N}\left[\exp(i \clv_N^{\varphi} + J^{N,1}(T) + J^{N,2}(T))\right]\\
 &\quad=\E\left[\exp(i I_1(\bar \Phi) + \frac{1}{2}I_2(F) - \frac{1}{2} \bar \tau)\right]\\
&\quad= \E \left[\E\left(\exp(i I_1(\bar \Phi) + \frac{1}{2}I_2(F) - \frac{1}{2} \bar \tau) \mid \clg^0_*\right)\right]\\
 &\quad=\int_{\Omega_m} \exp\left(-\frac{1}{2}(\sigma_{\om_0}^{\varphi})^2\right )P_0(d\omega_0),
 \end{align*}
 where the last equality is a consequence of Lemma 1.3 of \cite{ShiTan} and Lemma \ref{lem:lem916}.  Thus we have proved \eqref{eq:eq209} which completes the proof of Theorem
\ref{thm:clt}.
 \qed
\section{Convergence of the Signed Measures in the Path Space.}
\label{sec:pathspace}
In \cite{KuXi2} authors studied a functional central limit theorem for scaled and centered empirical measures for a family of weakly interacting particle systems with a common factor.
As noted in the Introduction, in the current work our focus is on limit theorems for functionals of the whole path of the particles, however in this section we will discuss how functional central limit 
theorems of the form in \cite{KuXi2} can be recovered from Theorem \ref{thm:clt}.
For $t \in [0,T]$ consider the random signed measure on $\RR^d$ defined as
\begin{equation}
	\label{eq:signmzr}
	\Lambda^N_t = \sqrt{N} \left( \frac{1}{N} \sum_{j=1}^N \delta_{Z^{N,j}_t} - \mu_t\right).
\end{equation}
We note that $\mu_t = \eta_t(\bar V^0)$ where $\bar V^0$ is as introduced in Section \ref{sec:cltassu} and for $\om_0 \in \Om_m$, $\eta_t(\om_0) = \alpha (\om_0, \cdot) \circ X_{*,t}^{-1}$ with
$\alpha$ as in \eqref{eq:transkern} and $X_*$ as in \eqref{eq:canon2}.

For notational simplicity we assume for rest of the section that $d=1$.
Following \cite{HitMot} and \cite{KuXi2} $\Lambda^N = \{\Lambda^N_t\}_{t\in[0,T]}$ can be regarded as a sequence of 
$\DD_{\Psi'}[0,T]$ valued random variables where $\Psi'$ is the dual of the 
``modified Schwartz space" $\Psi$  given as follows.
Let $\rho: \RR \to \RR$ be defined as
$$\rho(x) = C \exp \{ -1/(1-|x|^2)\} 1_{|x| < 1}, \; x \in \RR,$$
where $C \in (0,\infty)$ is such that $\int \rho(x) dx =1$.  Let
$$v(x) = \int e^{-|y|} \rho(x-y) dy, \; e(x) = 1/v(x), \; x \in \RR.$$
Let $\Psi = \{ \psi = e u: u \in \cls\}$ where $\cls$ is the Schwartz space (cf. \cite{GelVil}).  For $p \in \NN_0$ and $\psi \in \Psi$, define
$$
\|\psi\|_p^2 = \sum_{0\le k \le p} \int_R (1 + |x|^2)^{2k} \left| \frac{d^k}{dx^k}(\psi(x)v(x))\right|^2 dx.$$
Let $\Psi_p$ be the completion of $\Psi$ with respect to $\|\cdot\|_p$. The $\Psi_p$ is a Hilbert space with inner product $\langle \cdot, \cdot \rangle_p$
defined in an obvious manner.  For $\hat \phi \in \Psi_0$ and $\phi \in \Psi_{p}$
$$\hat \phi[\phi] \doteq \int_R \hat \phi(x) \phi(x) v^2(x) dx$$
defines a continuous linear functional on $\Psi_p$ with norm
$$\|\hat \phi\|_{-p} = \sup_{\phi \in \Psi_p} \frac{|\hat \phi[\phi]|}{\|\phi\|_p}.$$
Let $\Psi_{-p}$ be the completion of $\Psi_0$ with respect to this norm.  Then $\Psi$ is a nuclear space \cite{GelVil} and $\Psi' \doteq \cup_{k=0}^{\infty} \Psi_{-k}$ is its dual.

Recall the class $\cla$ introduced in Section \ref{sec:cltassu}. Given $\ell \in \NN$, $t_1, \cdots t_{\ell} \in [0,T]$ and $\phi_1, \cdots \phi_{\ell} \in \Psi$, define
$\bfvarphi_i \in \cla$, $i=1, \cdots, \ell$ as
$\bfvarphi_i(\om) = \phi_i(\om_{t_i})$, $\om \in \cld_d$.  Also, for $\om_0 \in \Om_m$, let $\Phi^i_{\om_0} = \bfvarphi_i(X_*) - m_{\bfvarphi_i}(\om_0)$
where $m_{\cdot}$ is as introduced in Section  \ref{sec:cltassu}.  Also define the $\ell \times \ell$ matrix
$\Sigma_{\om_0} = (\Sigma_{\om_0}^{ij})$, where
$$\Sigma_{\om_0}^{ij} = \langle (I-A_{\om_0})^{-1}\Phi_{\om_0}^i, (I-A_{\om_0})^{-1}\Phi_{\om_0}^j\rangle_{L^2(\Omega_d,\alpha(\omega_0,\cdot))}.$$
Let $\gamma^{\bfvarphi}_{\om_0}$ be a $\ell$ dimensional Gaussian random variable with mean $0$ and variance $\Sigma_{\om_0}$ and define
$$\gamma_{t_1, \cdots t_{\ell}}^{\phi_1, \cdots, \phi_{\ell}} \equiv \gamma ^{\bfvarphi} = \int_{\Om_m} \gamma^{\bfvarphi}_{\om_0} P_0(d\om_0).$$
The following theorem follows from Theorem \ref{thm:clt} of the current work and arguments similar to Theorem 3.1 of \cite{KuXi2}.  We only provide a sketch.
Let $\QQ^N \in \clp(D([0,T]:\Psi'))$ be the probability law of $\Lambda^N$. Define
$\Pi_{t_1, \cdots t_{\ell}}^{\phi_1, \cdots, \phi_{\ell}}: \DD_{\Psi'}[0,T] \to \RR^{\ell}$
as 
$$
\Pi_{t_1, \cdots t_{\ell}}^{\phi_1, \cdots, \phi_{\ell}}(u) = \left( u_{t_1}[\phi_1], \cdots, u_{t_{\ell}}[\phi_{\ell}]\right).
$$
\begin{theorem}
	\label{thm:pathcgce}
	Suppose all the assumptions in Theorem \ref{thm:clt} are satisfied. 
	Then, as $N \to \infty$, $\QQ^N \to \QQ$ where $\QQ$ is the unique probability measure on $\DD_{\Psi'}[0,T]$
	that satisfies
	$$\QQ \circ (\Pi_{t_1, \cdots t_{\ell}}^{\phi_1, \cdots, \phi_{\ell}})^{-1} = \gamma_{t_1, \cdots t_{\ell}}^{\phi_1, \cdots, \phi_{\ell}}$$
for all $\ell \ge 1$, $t_1, \cdots , t_{\ell}\in [0,T]$ and $\phi_1, \cdots \phi_{\ell} \in \Psi$.	
\end{theorem}
{\bf Sketch of Proof.}
From  Theorem 4.1 and Proposition 5.2 of \cite{Mit} it suffices to show that\\
\noindent 
(i) for every $\phi \in \Psi$, $\QQ^N_{\phi}$ is tight in 
$\DD_{\RR}[0,T]$, 
where $\QQ^N_{\phi} = \QQ^N\circ (\Pi^{\phi})^{-1}$ and $\Pi^{\phi}: \DD_{\Psi'}[0,T] \to \DD_{\RR}[0,T]$ is defined
as
$\Pi^{\phi}(u)[t] = u_t[\phi] $, $t \in [0,T]$.\\
\noindent
(ii) for all $\ell \ge 1$, $t_1, \cdots , t_{\ell}\in [0,T]$ and $\phi_1, \cdots \phi_{\ell} \in \Psi$,
$\QQ^N \circ (\Pi_{t_1, \cdots t_{\ell}}^{\phi_1, \cdots, \phi_{\ell}})^{-1} \to \gamma_{t_1, \cdots t_{\ell}}^{\phi_1, \cdots, \phi_{\ell}}$.

Proof of (i) follows along the lines of Theorem 3.1 of \cite{KuXi2} and is omitted.  Consider now (ii).  Fix $\ell \ge 1$, $t_1, \cdots , t_{\ell}\in [0,T]$ and $\phi_1, \cdots \phi_{\ell} \in \Psi$ as above.
Let $a_1, \cdots a_l \in \RR$ and define
$$\bfvarphi^{\bfa} = \sum_{i=1}^{\ell} a_i \bfvarphi_i, \; \Phi^{\bfa}_{\om_0} = \sum_{i=1}^{\ell} a_i \Phi^{i}_{\om_0},$$
where $\bfvarphi_i, \Phi^{i}_{\om_0}$ are as defined above the theorem.
Let $\sigma^{\bfa}_{\om_0} = \|(I-A_{\om_0})^{-1}\Phi_{\om_0}^{\bfa}\|_{L^2(\Omega_d,\alpha(\omega_0,\cdot))}$ and
$\tilde \gamma^{\bfa}_{\om_0}$ be a Normal random variable with mean $0$ and variance $(\sigma^{\bfa}_{\om_0})^2$ and let $\tilde \gamma^{\bfa} = \int \tilde \gamma^{\bfa}_{\om_0} P_0(d\om_0)$.
Let $\tilde \Pi^{\bfa}: \DD_{\Psi'}[0,T] \to \RR$ be defined as
$\tilde \Pi^{\bfa}(u) = \sum_{i=1}^{\ell} a_i u_{t_i}[\phi_i]$.
From Theorem \ref{thm:clt} it is immediate that $\QQ^N \circ (\tilde \Pi^{\bfa})^{-1} \to \tilde \gamma^{\bfa}$ as $N \to \infty$.
The statement in (ii) is now immediate from the classical Cram\'{e}r-Wold  argument.
\qed

\section{Application to Finance}
\label{sec:finapp}
Recently in \cite{CMZ}, authors have introduced a model for self-exciting correlated defaults in which default times of various entities depend not only on factors specific
to entities and a common factor but also on the average number of past defaults in the market.
The paper studies an asymptotic regime as the number of entities become large.  One of the results in \cite{CMZ} is a CLT which is established under somewhat
restrictive conditions on the model.  Below, we describe the result from \cite{CMZ} and then remark on how the results of current paper provide a CLT
for the model in \cite{CMZ} under much lesser restrictive conditions and for some of its variations.

\noindent The model for which CLT is considered in \cite{CMZ} (see Section 5.3 therein), using notation of the current paper, is as follows.
Let $(B^i)_{i\in \NN}$ be a sequence of  real standard Brownian motions and let $(\bfn^i)_{i\in \NN}$ be a sequence of Poisson random measures  on $\XX_T = [0,T] \times \mathbb{R} \times  \mathbb{R}_+$ with intensity measure $\bfnu = \lambda_T \otimes \delta_{\{1\}} \otimes \lambda_{\infty}$,
given on a   filtered probability space $(\Om, \clf, \PP, \{\clf_t\})$. All these processes are mutually independent and they have independent increments with respect to the filtration $\{\clf_t\}$.
  Consider the system of equations given by 
\begin{align}\label{finance}\begin{cases}
	U_t^{N}&=U_0 + \int_0^t\beta_0(U_s^N,\bar{L}^N_s)ds + \int_0^t\bar{\sigma}_0(U_s^N,\bar{L}^N_s)dB^0_s\,,\\
X_t^{N,i}&=X_0^{N,i} + \int_0^t\beta(X^{N,i}_s, Y^{N,i}_s, U_s^N,\bar{L}^N_s)ds + B^i_t\,,\,\,\,\,i=1,2\ldots, N\,,\\
Y_t^{N,i}&=\int_{\XX_t}1_{[0,\lambda(X^{N,i}_{s}, Y^{N,i}_{s-}, U_{s}^N,\bar{L}^N_{s-})]}(u)\bfn^i(ds\, dh\, du)\,,
\end{cases}
\end{align}
where $\bar{L}^N_t=\frac{1}{N}\displaystyle\sum_{i=1}^N\zeta(Y^{N,i}_t)$ for some bounded and Lipschitz map $\zeta$, and we assume that $\{X_0^{N,i}\}_{i=1}^N$ are i.i.d. with common distribution $\mu_0$ and $U_0$ is independent of $\{X_0^{N,i}\}_{i=1}^N$ and has probability distribution $\rho_0$.  Also, $\{X_0^{N,i}\}_{i=1}^N$ and $U_0$ are $\clf_0$ measurable. 
The  interpretation for the finance model is as follows. There are $N$ defaultable firms. The process $U^N$ represents the common factor process and $X^{N,i}$ is the $i$-th firm's  specific factor. $Y^{N,i}$ are counting processes representing
the number of defaults of firm $i$. The key feature of this model is that the correlation among the defaults not only depends on the common exogenous factor $U^N$, but also on the past defaults through the process $\bar{L}^N$. 
In the model of \cite{CMZ}, $\zeta(y) = |y| \wedge 1$ and consequently all values of $Y^{N,i}_t$ greater than $0$ are treated the same way (an entity has either not defaulted by time $t$ or it has  defaulted  in which case it disappears from the system.)  
The paper \cite{CMZ} establishes a CLT for $\bar{L}^N_t$ under the condition that
$\lambda(x, y, u, l) \equiv \lambda(l)$, $x, y, u, l \in \RR$.  Note that in this case the factor processes $X^{N,i}$  and $U^N$ become irrelevant.

The model in \eqref{finance} is a special case of the model considered in  \eqref{eq:nparteq} and \eqref{eq:comminp} with the following  identifications:
\begin{itemize}
\item $d=2$, $m=1$.
\item $Z^{N,i} = (X^{N,i}, Y^{N,i})'$.
\item $b=(\bar\beta,0)',\,\,b_0=\bar\beta_0$, where for $z \in\RR^2$, $u \in \RR$, $\nu \in \clp(\RR^2)$,
$$\bar \beta(z,u,\nu) = \beta(z,u, \langle  \hat \zeta, \nu \rangle), \; \bar \beta_0(u,\nu) = \beta_0(u, 
\langle  \hat \zeta, \nu \rangle),$$
where $\hat \zeta: \RR^2 \to \RR$ is defined as $\hat \zeta(x,y) = \zeta(y)$, $(x,y) \in \RR^2$.
\item $\sigma_0(u,\nu)=\bar{\sigma}_0(u, 
\langle  \hat \zeta, \nu \rangle)$, $\sigma=\biggl(\begin{array}{cc}
1 & 0 \\ 
0 & 0
\end{array}\biggr)$ (See Remark \ref{rem2.3mod}).
\item $d_0=0$, $d=\bar \lambda$, where for $z \in\RR^2$, $u \in \RR$, $\nu \in \clp(\RR^2)$,
$\bar \lambda(z,u,\nu) = \lambda(z,u, \langle \hat \zeta , \nu \rangle)$.
\end{itemize}
Coefficients $\beta, \beta_0$ and $\lambda$ are required to satisfy the following conditions.

\noindent\textbf{(A1)} The function $\beta$ is bounded and Lipschitz.
 $\beta(z,u,l)$ is  twice continuously differentiable in $u$ and $l$ with bounded derivatives. 

\noindent\textbf{(A2)} The function $\beta_0$ is bounded and Lipschitz.
 $\beta_0(u,l)$ is  twice continuously differentiable in $u$ and $l$ with bounded derivatives. Exactly same assumptions for $\bar{\sigma}_0$
 
\noindent\textbf{(A3)} The function $\lambda$ is nonnegative, bounded, Lipschitz and  it is bounded away from $0$. 
 $\lambda(z,u,l)$ is  twice continuously differentiable in $u$ and $l$ with bounded derivatives.

Under the above assumptions it can be easily checked that Conditions \ref{cond:main}, \ref{cond:clt2}, \ref{cond:clt3}, \ref{cond:clt4}  and the modified form of Condition \ref{cond:clt1} in Remark \ref{rem2.3mod} are satisfied. Thus from Theorem \ref{thm:clt} it follows
that the average default process $\{\bar{L}^N_t\}$ satisfies a CLT. More precisely, for $t \in [0,T]$,
$$\sqrt{N}(\bar{L}^N_t-m_t(B^0,U_0))$$ converges in distribution to a random variable whose distribution is given as a mixture of Gaussians, where for $t \in [0,T]$,
$m_t: \clc_1 \times \mathbb{R} \to [0,1]$ is the measurable map such that $m_t(B^0, U_0) = \EE(\zeta(Y_t) \mid B^0, U_0)$ if
$(U, X, Y, \alpha)$ solve the following nonlinear system of equations.
%
%
%
%
%
\begin{align*}U_t&=U_0 + \int_0^t\beta_0(U_s,\alpha_s)ds + B^0_t\,,\\
X_t&=X_0 + \int_0^t\beta(X_s,Y_s, U_s,\alpha_s)ds + B_t\,,\\
Y_t&=\int_{\XX_t}1_{[0,\lambda(X_{s-},Y_{s-}, U_{s-},\alpha_{s-})]}(u)dN,\; \alpha_t = \EE(\zeta(Y_t) \mid B^0, U_0),
\end{align*}
where $B^0$ and $B$ are  Brownian motions and $N$ is a Poisson random measure on $\XX_T$ with intensity measure $\bfnu$,
given on  $(\Om, \clf, \PP, \{\clf_t\})$ such that they are mutually independent and they have independent increments with respect to the filtration $\{\clf_t\}$. Also,  $X_0$ and $U_0$ are independent $\clf_0$ measurable random variables with distribution $\mu_0$  and  $\rho_0$ respectively.  

The results of the current paper (in contrast to \cite{CMZ}) not only allow for a general dependence of $\lambda$
on factor processes but can also be used to treat more complex forms of default processes and also settings where the common factor and specific factor dynamics have both diffusion and jump components.  




\bibliographystyle{plain}

{\sc
\bigskip
\noi
A. Budhiraja\\
Department of Statistics and Operations Research\\
University of North Carolina\\
Chapel Hill, NC 27599, USA\\
email: budhiraj@email.unc.edu
\skp

\noi
E. Kira\\
Department of Statistics\\
Institute of Mathematics and Statistics\\
University of Sao Paulo\\
Rua do Matao 1010, 05508-090, Sao Paulo, SP, Brazil\\
email: betikira@ime.usp.br 

\skp

\noi
S.Saha\\
Department of Electrical Engineering\\
Technion - Israel Institute of Technology\\
Haifa 32000, Israel\\
email: subhamay@tx.technion.ac.il

}

\end{document}